\documentclass[12pt,reqno, article]{amsart}
\usepackage{setspace}
\usepackage[foot]{amsaddr}
\usepackage{bm}
\usepackage{tikz-qtree}
\usepackage{booktabs}
\usepackage{multirow}
\usepackage{bbm}
\usepackage{diagbox}
\usetikzlibrary{positioning}
\usepackage[per=slash]{siunitx}
\usepackage[plain,noend]{algorithm2e}

\usepackage[left = 1in, right=1in, top = 1in, bottom = 1 in]{geometry}

\newtheorem{theorem}{Theorem}
\newtheorem{proposition}{Proposition}
\newtheorem{corollary}{Corollary}
\newtheorem{lemma}{Lemma}

\theoremstyle{definition}
\newtheorem{condition}{Condition}

\theoremstyle{remark}

\newtheorem{remark}{Remark}
\renewcommand{\P}{{\text{pr}}}

\newcommand{\1}{1}

\newcommand{\T}{{\mathtt{T}}}

\newcommand{\cP}{{\mathcal{P}}}

\def\spacingset#1{\renewcommand{\baselinestretch}%
{#1}\small\normalsize} \spacingset{1}

\usepackage{natbib}
\bibpunct{(}{)}{;}{a}{}{,}

\allowdisplaybreaks

\title[]{Prepivoted permutation tests}
\author{Colin B. Fogarty$^{1}$} \address{$^{1}$ Operations Research and Statistics Group, MIT Sloan School of Management, Massachusetts Institute of Technology, Cambridge, MA }

\begin{document}

\begin{abstract}
We present a general approach to constructing permutation tests that are both exact for the null hypothesis of equality of distributions and asymptotically correct for testing equality of parameters of distributions while allowing the distributions themselves to differ. These robust permutation tests transform a given test statistic by a consistent estimator of its limiting distribution function before enumerating its permutation distribution. This transformation, known as prepivoting, aligns the unconditional limiting distribution for the test statistic with the probability limit of its permutation distribution. Through prepivoting, the tests permute one minus an asymptotically valid $p$-value for testing the null of equality of parameters. We describe two approaches for prepivoting within permutation tests, one directly using asymptotic normality and the other using the bootstrap. We further illustrate that permutation tests using bootstrap prepivoting can provide improvements to the order of the error in rejection probability relative to competing transformations when testing equality of parameters, while maintaining exactness under equality of distributions. Simulation studies highlight the versatility of the proposal, illustrating the restoration of asymptotic validity to a wide range of permutation tests conducted when only the parameters of distributions are equal.
 
\end{abstract}




\maketitle
\thispagestyle{empty}

\spacingset{1.8} 

\setcounter{page}{1}
\section{Introduction}
Permutation tests are commonly used for inference in applications ranging from field trials in economics to microarray analyses in genomics. In many applied domains, permutation tests are erroneously viewed as appealing nonparametric substitutions for the usual large-sample tests assessing equality of parameters of distributions \citep{lud98, rux06}.  Unfortunately, while permutation tests provide exact tests for the null hypothesis of equality of distributions, they do not generally provide valid inference for the equality of parameters of distributions even asymptotically. For instance, a permutation test which permutes the difference in sample means does not generally provide an asymptotically valid test for equality of population means \citep{rom90}.  While it may be correct avoid the $t$ test in small sample regimes, replacing it with a permutation test generally requires a reformulation of the null hypothesis being tested. Absent this, practitioners may draw unwarranted conclusions from a permutation test's rejection of the null hypothesis. 


Misconceptions surrounding the findings to which practitioners are entitled motivate the development of robust permutation tests. Here robustness does not refer to any particular connection with the field of robust statistics, but rather to a desired robustness against natural misinterpretations when implementing a permutation test and assessing its results. A robust permutation test maintains the exactness of the procedure under the null of equality in distributions while, through a suitably chosen test statistic, providing asymptotically valid inference for the null hypothesis that certain parameters of the distributions are equal even when the distributions differ. Building upon work of \citet{neu93}, \citet{jan97}, \citet{pau11}, \citet{ome12} among others, \citet{chu13} showed in the univariate setting that when asymptotically linear estimators exist for a parameter $\theta(\cdot)$, simply permuting a studentized version of these estimators suffices to yield an asymptotically robust permutation test. \citet{chu16mult} extend these results to the multivariate setting, suggesting permuting a modified Hotelling T-squared statistic using an unpooled covariance estimator. The success of studentization in the aforementioned problems is tied to the standard normal (univariate) and chi-squared (multivariate) limiting distributions attained after studentization. Both the unconditional asymptotic distribution of the test statistic and the limit of the reference distribution generated by permutation align after studentization, resulting in asymptotic correctness of permutation inference for testing equality of parameters even when the distributions themselves differ. Studentization only produces a robust permutation test when it yields an asymptotic pivot for both the unconditional and permutation distributions, hence failing to provide a general resolution. 

Prepivoting is the transformation of a test statistic by an empirical estimate of its cumulative distribution function. The transformation is motivated by the probability integral transform: if we had access to the true distribution function $H_{T_{n}}(x, P)$ of a continuous random variable $T_{n}$, then $H_{T_{n}}(T_{n}, P)$ would be distributed as a uniform random variable on [0,1]. We describe two forms of prepivoting that can restore asymptotic validity to permutation tests when testing the null of equality for parameters based upon estimators admitting asymptotically linear representations. The first, termed Gaussian prepivoting, makes direct use of the asymptotic normality stemming from asymptotic linearity. It is shown that this transformation is equivalent to studentization whenever studentization would produce an asymptotic pivot, but provides an asymptotic pivot for a broader class of test statistics. Gaussian prepivoting both sheds light upon the true purpose of studentization when it is successful and provides general methodology for situations where studentization alone would fail. \citet{coh21gauss} describe Gaussian prepivoting within randomization tests in finite population causal inference as a way to create tests that are exact for Fisher's sharp null while asymptotically conservative for Neyman's weak null. The second, bootstrap prepivoting, forms an estimate of the distribution function using the bootstrap. Bootstrap prepivoting was introduced in \citet{ber87,ber88} as a means of improving the order of approximation for bootstrap confidence intervals and hypothesis tests. \citet{chu16mult} first used bootstrap prepivoting within a permutation test to restore first-order correctness to multivariate permutation tests permuting the maximum absolute $t$ statistic. In this work we highlight the broad generality with which the idea succeeds, along with circumstances in which bootstrap prepivoting provides higher order accuracy relative to competing methods for constructing robust permutation tests. 


This work suggests the following general approach to constructing robust permutation tests: rather than permuting a seemingly natural test statistic for assessing equality of parameters, instead permute one minus a large-sample $p$-value known to be both asymptotically valid under the null of equality of parameters and consistent in level should the distributions themselves be equal. This straightforward prescription has appeal even for researchers solely interested in testing equality of parameters who already have access to an asymptotically valid $p$-value for said task. By basing inference on the permutation distribution of one minus the $p$-value  rather than the $p$-value itself as is convention, one guarantees exactness of the resulting inference at any sample size under equality of distribution without degrading performance under equality of parameters. The cost of conferring this exactness property lies not in the development of a clever test statistic (one minus an existing large-sample $p$-value-is permuted), but rather in additional computation time from conducting the permutation test. Theoretical developments and simulation studies in \S\S\ref{sec:higher}-\ref{sec:sim} highlight the particular benefits that bootstrap $p$-values may confer in terms of reduced errors in rejection probability.

\section{Notation and review}
\subsection{Constructing a permutation distribution}
Let $P_1,\ldots,P_k$ be ${d}_P$-variate distributions for $k$ populations. For each $i=1,\dots, k$ let $X_{i1},\ldots,X_{in_i}$ be independent and identically distributed random variables distributed according to $P_i$, and let the $k$ samples be jointly independent. Let $n = \sum_{i=1}^kn_i$, and impose the following condition on $n_i/n$:
\begin{condition}\label{cond:frac} For $i=1,...,k$, $n_i/n\rightarrow p_i$ with $0 < p_i < 1$ as $n\rightarrow \infty$, with $n_i/n - p_i = O(n^{-1/2})$.
\end{condition} Let $X_i$ denote the $n_i\times {d}_P$ matrix whose $j$th row contains $X_{ij}$. Write $Z$ for the $n\times {d}_P$ matrix stacking $X_1,\ldots,X_k$ on top of one another, and let $\mathcal{I}_i$ be the row indices of $Z$ containing samples from population $i$, i.e. $\mathcal{I}_1 = \{1,\ldots,n_1\},  \mathcal{I}_2 = \{n_1+1,\ldots,n_1+n_2\}$,\ldots, $\mathcal{I}_k = \{n - n_k+1,\ldots,n\}$. Let $P^n = \prod_{i=1}^k\prod_{j=1}^{n_i}P_i$ be the distribution of $Z$, and let $P^\infty$ be the corresponding infinite product measure sending $n\rightarrow \infty$. We will generally drop the dependence on the sample size and refer to both $P^n$ and $P^\infty$ as $P$. 

Let $\pi = \{\pi(1),\ldots,\pi(n)\}$ be a permutation of $\{1,\dots,n\}$. Let $Z_\pi$ reflect the matrix $Z$ with its rows permuted according to $\pi$, so that the $j$th row of $Z_\pi$ contains the $\pi(j)$th row of $Z$. Let $Z_{\pi i} = Z_{\pi(\mathcal{I}_i)}$ be the $n_i\times {d}_P$ matrix containing the $n_i$ rows of $Z_\pi$ whose row indices are within $\mathcal{I}_i$. For instance, the $j$th row of $Z_{\pi i}$ contains the $\pi\{\mathcal{I}_i(j)\}$th row of $Z$. Define $\hat{P}_i$ as the empirical distribution of $Z_i$, $\hat{P} = \prod_{i=1}^k\prod_{j=1}^{n_i}\hat{P}_{i}$, and let the analogous definitions hold for $\hat{P}_{\pi i}$ and $\hat{P}_\pi$ applied to $Z_\pi$. The permutation distribution for a statistic $T_{n}(Z_\pi)$ is
\begin{align*}
R_{T_{n}}(x, \hat{P}) &= \frac{1}{n!}\sum_{\pi \in \cP_{n}}\1\{T_{n}(Z_{\pi})\leq x\} = \P\{T_{n}(Z_\Pi)\leq x\mid Z\},
\end{align*} where $\mathcal{P}_n$ is the set of all $n!$ permutations $\pi$, $\Pi$ is uniformly distributed on $\cP_{n}$ and is independent of $Z$, and $\1(A)$ is an indicator that the event $A$ occurred. 

Consider the null hypothesis of equality of distributions $\mathcal{H}_P: P_1=\ldots=P_k$. The level$-\alpha$ permutation test for $\mathcal{H}_P$ based upon a scalar test statistic $T_n$ rejects the null hypothesis through the decision rule \begin{align*}\varphi_{T_{n}}(\alpha) &= \1\{T_{n}(Z) > R_{T_{n}}^{-1}(1-\alpha,\hat{P})\},
\end{align*} where $R_{T_n}^{-1}(1-\alpha, \hat{P}) = \inf\{x: R_n(x, \hat{P}) \geq 1-\alpha\}$. Through this construction, $E\{\varphi_{T_{n}}(\alpha)\} \leq \alpha$ for any $n$. Control in finite samples stems from the order statistics forming a complete sufficient statistic for the common unknown distribution under $\mathcal{H}_P$. It follows that only exact tests have a Neyman structure, and the distribution of $T_{n}(Z)$ given the order statistics is precisely the permutation distribution $R_{T_{n}}(\cdot, \hat{P})$ \citep[Theorem 5.8.1]{leh05}. A level exactly equal to $\alpha$ may be attained through a suitably randomized decision rule when $T_{n}(Z) = R_n^{-1}(1-\alpha, \hat{P})$. 
 
\subsection{Testing equality of parameters}
 Suppose interest lies not in equality of distributions themselves but rather in the null of equality of parameters of distributions. Let $\theta(\cdot) = (\theta_1(\cdot),...,\theta_d(\cdot))^\T$ be a $d$-dimensional parameter and consider testing the null  \begin{align*} \mathcal{H}_{\theta(P)}: \theta_\ell(P_1) = \cdots = \theta_\ell(P_k)\;\; \text{for all } \ell = 1,\ldots,d,\end{align*}against the alternative that $\theta_\ell(P_i)\neq \theta_\ell(P_i')$ for some $\ell = 1,\ldots,d$ and some $i\neq i'$. Note that $d$ may differ from ${d}_P$, the dimension of the distributions $P_i$. We suppose the existence of estimators $\hat{\theta}_\ell(\cdot)$ for each $\ell$ which satisfy the following:
 \begin{condition}\label{cond:asymlin}
There exist asymptotically linear estimators $\hat{\theta}_\ell(\cdot)$ such that for each $i=1,\ldots,k$ and each $\ell=1,\ldots,d$
\begin{align*}
n_i^{1/2}\{\hat{\theta}_\ell(X_{i}) - \theta_\ell(P_i)\} &= n_i^{-1/2}\sum^{n_i}_{j=1}f_{P_{i}, \ell}(X_{ij}) + o_{P_i}(1),
\end{align*} 
with $E\{f_{P_{i},\ell}(X_{i1}))\} = 0$ and $E\{f_{P_{i}, \ell}(X_{i1})^2\} < \infty$.  Furthermore, let $Y_1,\ldots,Y_n$ be independent and identically distributed from the mixture $\bar{P} = \sum_{i=1}^kp_{i}P_i$ and let $Y$ be the $n\times {d}_P$ matrix whose $j$th row is $Y_{j}$. Then for each $\ell= 1,...,d$,
\begin{align*}
n^{1/2}\{\hat{\theta}_\ell(Y) - \theta_\ell(\bar{P})\} &= n^{-1/2}\sum^{n}_{j=1}f_{\bar{P}, \ell}(Y_{j}) + o_{\bar{P}}(1),
\end{align*} with $E\{f_{\bar{P},\ell}(Y_{1}))\} = 0$ and $E\{f_{\bar{P}, \ell}(Y_{1})^2\} < \infty$.
\end{condition}
Let $\Sigma_i$ be the covariance matrix for $(f_{P_{i},1}(X_{i1}),\ldots, f_{P_{i},d}(X_{i1}))^\T$ for $i=1,..,k$, and let $\bar{\Sigma}$ be the covariance matrix for $(f_{\bar{P},1}(Y_{1}),\ldots, (f_{\bar{P},d}(Y_{1}))^\T$ with $Y_1$ distributed as in Condition \ref{cond:asymlin}.  Let $\hat{\Theta}(Z_\pi)$ be the $d\times k$ matrix with $\hat{\theta}_\ell(Z_{\pi i})$ in its $\{\ell,i\}$ entry. We will generally use $\hat{\Theta} = \hat{\Theta}(Z)$ and $\hat{\Theta}_\pi = \hat{\Theta}(Z_\pi)$ as shorthand  for these and other estimators.

Let $\hat{C}_\pi$ be a matrix of column contrasts of dimension $k \times m$, i.e. a matrix whose $m$ columns all sum to zero with no column equal to the zero vector. We begin with test statistics $T_{n}$ satisfying the following:
\begin{condition}\label{cond:T}
$T_n(\cdot)$ is of the form
\begin{align}\label{eq:T}
T_n(Z_\pi) &= g\{n^{1/2}\hat{\Theta}_\pi\hat{C}_\pi, \hat{\eta}_\pi\},
\end{align}
where (a) $g: \mathbb{R}^{d\times m} \times \Xi \mapsto \mathbb{R}$ is a jointly continuous function in both of its arguments; (b) $\hat{C}$ and $\hat{C}_\Pi$ converge in probability to contrast matrices $C$ and $\bar{C}$ respectively; and (c) $\hat{\eta}$ and $\hat{\eta}_\Pi$ converge in probability to values $\eta, \bar{\eta} \in \Xi$ respectively, where $\Pi$ is uniformly distributed over $\mathcal{P}_n$ and is independent of $Z$.
\end{condition}

As will be demonstrated in \S\ref{sec:sim}, many common test statistics for testing equality in parameters may be written in this form. Without loss of generality, in what follows we assume that we reject the null for large values of $T_n$. Let $H_{T_{n}}(x, P)$ be the true distribution function for $T_{n}$ and suppose that $\mathcal{H}_{\theta(P)}$ holds but $P_i\neq P_{i'}$ for some $i,i' = 1,\ldots,k$. Let $H_{T}(x, P)$ be the limiting distribution function for $T_{n}$, and suppose that the permutation distribution $R_{T_{n}}(x, \hat{P})$ converges in probability to $R_T(x, P)$ at all continuity points of $R_T(x, P)$. Then, even for seemingly natural test statistics $T_n$ satisfying Conditions \ref{cond:asymlin} and \ref{cond:T}, it can occur that $\sup_{x\in\mathbb{R}}|H_T(x, P)- R_T(x, P)| > 0$. As a result, the permutation test $\varphi_{T_{n}}(\alpha)$ may be anti-conservative for $\mathcal{H}_{\theta(P)}$ even in the limit. Perhaps the most prominent example of this is the permutation distribution of the difference in means in the two-sample univariate case, which fails to control the Type I error rate even asymptotically under the null of equality of expectations of distributions unless either the variances of the two distributions are equal or the sample sizes are equal; see \citet{rom90} or \citet[Example 2.1]{chu13}.


This potential for anti-conservativeness for statistics satisfying Conditions \ref{cond:asymlin} and \ref{cond:T} may be understood through the following theorem comparing the unconditional  limit distribution of $n^{1/2}\text{vec}(\hat{\Theta}\hat{C})$ to the probability limit of its permutation distribution, where $\text{vec}(\cdot)$ is the columnwise vectorization operator.

\begin{theorem}\label{thm:CLT}
Suppose Conditions \ref{cond:frac} and \ref{cond:asymlin} hold and let $V_n(Z_\pi) = n^{1/2}\text{vec}(\hat{\Theta}_\pi\hat{C}_\pi)$.  Unconditionally, $V_n(Z)$ converges in distribution to a multivariate normal random variable, with mean zero and covariance 
\begin{align}\label{eq:cov}
(C^\T\otimes I_{d\times d})\Gamma(C^\T \otimes I_{d\times d})^\T,
\end{align}
where $\otimes$ is the Kronecker product, $I_{d\times d}$ is the $d\times d$ identity matrix, and $\Gamma$ is the direct sum of the matrices $p_i^{-1}\Sigma_i$ $(i=1,..,k)$, i.e. the block-diagonal matrix of dimension $kd\times kd$ with $p_i^{-1}\Sigma_i$ in the $i$th block of dimension $d\times d$. 

The permutation distribution of $V_n(\cdot)$, $R_{V_{n}}(t,  \hat{P})$, instead converges weakly in probability to the law of a multivariate normal distribution, with mean zero and covariance 
\begin{align}\label{eq:covperm}
(\bar{C}^\T\otimes I_{d\times d})\bar{\Gamma}(\bar{C}^\T\otimes I_{d\times d})^\T,
\end{align}
where $\bar{\Gamma}$ is the direct sum of the matrices $p_i^{-1}\bar{\Sigma}$ $(i=1,\ldots,k)$.
\end{theorem}

Under  $\mathcal{H}_{\theta(P)}$ $\bar{\Gamma}$ and $\bar{C}$ need not equal $\Gamma$ and $C$, such that the covariances governing the true limit distribution and the probability limit of the permutation distribution generally differ. As the true limiting distribution and the limiting permutation distribution for $T_n$ are simply push-forward measures of the above multivariate normals determined by the function $g$ in (\ref{eq:T}), this can corrupt the limiting size of $\varphi_{T_n}$ under $\mathcal{H}_{\theta(P)}$. That is, Condition \ref{cond:T} is insufficient for the permutation distribution of a test statistic $T_n$ to provide asymptotically valid inference for $\mathcal{H}_{\theta(P)}$.
\begin{remark}
Given asymptotic linearity, establishing asymptotic normality for the unconditional distribution of $V_n(Z)$ is a standard application of the central limit theorem and is omitted. The proof of weak convergence in probability of $R_{V_n}$ to the distribution function of the mean zero multivariate Gaussian with covariance given in (\ref{eq:covperm}) is an extension of steps within the proof of Theorem 3.1 of \citet{chu13} to the multivariate setting. In sketching the proof in the web-based appendix, we highlight the essential role of the contrast matrix $\hat{C}(Z_\pi)$ in the definitions of $V_n(Z_\pi)$ and $T_n(Z_\pi)$. Under our assumptions, $R_{V_n}$ converges weakly in probability to the law of a mean-zero multivariate normal with covariance given above if and only if $\hat{C}_n(Z_\Pi)$ converges in probability to a contrast matrix. In short, the contrast matrix justifies application of Lemma 4.1 of \citet{dum13}, as the independence required in their (D2) is satisfied; see also Lemma A.1 of \citet{chu16mult}.  Example 5.3 of \citet{chu13} shows what occurs when $\hat{C}_\Pi$ does not converge in probability to a contrast matrix: outside the degenerate case with all columns equal to zero, the permutation distribution will not converge weakly in probability to any fixed law $Q$.

\end{remark}

\section{Prepivoted permutation tests}
\subsection{A general approach through the probability integral transform}
Let $H_{T}(x,P)$ be the limiting law for a scalar test statistic $T_{n}(Z)$, let $R_T(x, P)$ represent the probability limit of $R_{T_{n}}(x,\hat{P})$, and suppose that both $H_T(x,P)$ and $R_T(x,P)$ are continuous and strictly increasing. Suppose that $F_{T_{n}}(x,  \hat{P})$ is an estimate of $H_{T_{n}}(x,P)$, the distribution function for $T_{n}(Z)$. Consider now the true distribution $H_{F_{n}}$ of $F_{T_{n}}\{T_{n}(Z),  \hat{P}\}$, along with the permutation distribution $R_{F_{n}}$ of the form\begin{align*}
R_{F_{n}}(x,  \hat{P}) &= \frac{1}{n!}\sum_{\pi \in \cP_{n}}\1[F_{T_{n}}\{T_{n}(Z_{\pi}), \hat{P}_{\pi}\}\leq x] =  \P[F_{T_{n}}\{T_{n}(Z_\Pi),\hat{P}_\Pi\} \leq x \mid Z].
\end{align*} In words, $F_{T_{n}}\{T_{n}(Z),  \hat{P}\}$ transforms $T_{n}(Z)$ by an estimate of its distribution function, in so doing approximating a transformation that would  give rise to a random variable distributed uniformly on [0,1] under continuity. Theorem \ref{thm:general} shows that if $F_{T_n}(x, \hat{P}_\Pi)$ also suitably estimates $R_T(x,P)$ the limiting distributions $H_{F}$  and $R_F$ will be the same.

\begin{theorem}\label{thm:general}
Let $\mathcal{A}_{\theta(P)}$ be a subset of the distributions $P_1,...,P_k$ satisfying $\mathcal{H}_{\theta(P)}$. Suppose that for each $(P_1,\ldots,P_k) \in \mathcal{A}_{\theta(P)}$ the limits $R_T(x,P)$ and $H_T(x,P)$ are continuous and strictly increasing. Suppose  that the estimator $F_{T_{n}}$ satisfies 
\begin{align*}
\underset{x \in \mathbb{R}}{\sup} |F_{T_{n}}(x,  \hat{P}) - H_{T_{n}}(x, P)| \rightarrow 0;\;\;\;\underset{x \in \mathbb{R}}{\sup} |F_{T_{n}}(x,\hat{P}_\Pi) - R_{T_{n}}(x,\hat{P})| \rightarrow 0\end{align*} in probability as $n\rightarrow\infty$ for each element of $\mathcal{A}_{\theta(P)}$. Then, for each element of $\mathcal{A}_{\theta(P)}$\begin{align*}
\underset{x \in \mathbb{R}}{\sup} |H_{F_{n}}(x, P) - U(x)| \rightarrow 0;\;\;\;\underset{x \in \mathbb{R}}{\sup} |R_{F_{n}}(x,\hat{P}) - U(x)| \rightarrow 0 \end{align*} in probability, where $U(x)$ is the distribution function for a uniform random variable on [0,1]. 
\end{theorem}

The theorem indicates that after transformation by $F_{T_{n}}$, the true distribution $H_{F_{n}}$ and the permutation distribution $R_{F_n}$ both converge to the same limits $H_F(x,P) = R_F(x,P) = U(x)$, even if the limits $H_{T}$ and $R_T$ are not the same. Applying Lemma 11.2.1 and Corollary 11.2.3 of \citet{leh05} yields the following corollary:
\begin{corollary}
Suppose the setup of Theorem \ref{thm:general} and consider the permutation test 
\begin{align*}
\varphi_{F_n}(\alpha) = \1[F_n\{T_{n}(Z),\hat{P}\} > R_{F_n}^{-1}(1-\alpha,\hat{P})].
\end{align*}
Then for each $(P_1,\ldots,P_k) \in \mathcal{A}_{\theta(P)}$, $\lim\;E\{\varphi_{F_n}(\alpha)\} = \alpha$. That is, the test $\varphi_{F_n}(\alpha)$ is asymptotically pointwise level $\alpha$ for all elements of $\mathcal{A}_{\theta(P)}$.
\end{corollary}

\subsection{Permuting asymptotically valid $p$-values}
Suppose inference is desired for $\mathcal{H}_{\theta(P)}$  but that the requirement of exactness when $P_1=\ldots=P_k$ is dropped. If $F_{T_{n}}$ satisfies $\sup_{x\in\mathbb{R}}\;\; |F_{T_{n}}(x,  \hat{P}) - H_{T_{n}}(x, P)| \rightarrow 0$ in probability for any element of  $\mathcal{A}_{\theta(P)}$ with a continuous and strictly increasing limiting distribution function $H_T(x, P)$, then an asymptotically valid $p$-value over this subset of distributions satisfying $\mathcal{H}_{\theta(P)}$ is
\begin{align*}
p_{val} &= 1-F_{T_{n}}\{T_{n}(Z),\hat{P}\}.
\end{align*}
Rejecting the null if $1-F_{T_{n}}\{T_{n}(Z),  \hat{P}\} \leq \alpha \Leftrightarrow F_{T_{n}}\{T_{n}(Z),  \hat{P}\}\geq 1-\alpha$ provides a pointwise asymptotically level $\alpha$ test over $\mathcal{A}_{\theta(P)}$. The permutation test $\varphi_{F_n}(\alpha)$ instead rejects when
\begin{align*}F_{T_{n}}\{T_{n}(Z), \hat{P}\} > R^{-1}_{F_{n}}(1-\alpha,  \hat{P}).
\end{align*}
The permutation distribution $R_{F_{n}}(x,  \hat{P})$ calculates for each permutation $\pi$ what 1 minus the large-sample $p$-value using $F_{T_{n}}$ would be with an observed sample $Z_{\pi}$. If ${\sup}_{x\in\mathbb{R}} |F_{T_{n}}(x,\hat{P}_\Pi) - R_{T_{n}}(x,\hat{P})| \rightarrow 0$ in probability, Theorem \ref{thm:general} and Lemma 11.2.1 of \citet{leh05} imply that $R^{-1}_{F_{n}}(1-\alpha,  \hat{P})\rightarrow U^{-1}(1-\alpha) = 1-\alpha$ in probability. In words, the $1-\alpha$ quantile of the permutation distribution for one minus the large-sample $p$-value converges in probability to $1-\alpha$. While the usual test would reject when one minus the $p$-value constructed using $F_{T_{n}}$ falls at or above $1-\alpha$, the permutation test $\varphi_{F_n}(\alpha)$ replaces $1-\alpha$ with the permutation quantile $R^{-1}_{F_{n}}(1-\alpha,\hat{P})$, which converges in probability to $1-\alpha$. Replacing $1-\alpha$ by the $1-\alpha$ quantile of this permutation distribution preserves asymptotic correctness under $\mathcal{A}_{\theta(P)} \subset \mathcal{H}_{\theta(P)}$, while additionally providing exactness should $P_1=\ldots=P_k$. 

Given an asymptotically valid $p$-value, constructing a robust permutation test under a given restriction on $\mathcal{H}_{\theta(P)}$ amounts to confirming that the permutation distribution of the large-sample $p$-value converges weakly in probability to the distribution function of a uniform random variable on [0,1]. This task is greatly simplified through the theoretical advances presented in \citet{chu13} and \citet{chu16mult}. While many approaches exist for constructing robust permutation tests that are tailored to specific test statistics $T_{n}$, we now present two general approaches for constructing asymptotically valid $p$-values whose permutation distributions have the required limiting behavior. The first approach, presented in \S\ref{sec:asym}, makes direct use the asymptotic normality of $V_n(Z)$ and $V_n(Z_\Pi)$ proven in Theorem \ref{thm:CLT} to compute an asymptotically valid $p$-value using a suitably constructed covariance estimator. The second, presented in \S \ref{sec:boot}, instead permutes one minus a large-sample $p$-value using a bootstrap null distribution for $\mathcal{H}_{\theta(P)}$ . 


\subsection{Gaussian prepivoting}\label{sec:asym}
Our first approach uses the asymptotic normality proven in Theorem \ref{thm:CLT}. We require the existence of a covariance estimator $\hat{\Sigma}$ of the following form:
\begin{condition}\label{cond:sigma} For each $i$, $\hat{\Sigma}(X_i)$ converges in probability to $\Sigma_i$. Furthermore, $\hat{\Sigma}(Y)$ converges in probability to $\bar{\Sigma}$, where $Y_1,..,Y_n$ are independent and identically distributed from the mixture $\bar{P} = \sum_{i=1}^kp_iP_i$ and $Y$ is an $n\times {d}_P$ matrix whose $j$th row is $Y_j$.  
\end{condition} For any permutation $\pi \in \mathcal{P}_n$, let $\hat{\Gamma}_\pi = \hat{\Gamma}(Z_\pi)$ be the block-diagonal matrix of dimension $kd\times kd$ with $\hat{\Sigma}(Z_{\pi i})$ in the $i$th  of $k$ blocks. Consider the distributional estimator
\begin{align}\label{eq:asym}K_{T_n}(x, \hat{P}_\pi) &= \gamma^{(kd)}_{0, \hat{\Gamma}_\pi}\left\{{a}: {g}(\text{vec}_{d,k}^{-1}(a)\hat{C}_\pi, \hat{\eta}_\pi) \leq x \right\},\end{align} where $\text{vec}_{d,k}^{-1}(a)$ is a map from $(dk) \times 1$ vectors $a$ to $d \times k$ matrices $A$ defined by $\text{vec}_{d,k}^{-1}\{\text{vec}(A)\} =A$ and $\text{vec}\{\text{vec}_{d,k}^{-1}(a)\}=a$, where $a= \text{vec}(A)$; and $\gamma^{(p)}_{\mu, \Lambda}(\mathcal{A})$ is the probability that a $p$-variate Gaussian random variable with mean parameter $\mu$ and covariance parameter $\Lambda$ falls within a set $\mathcal{A}$, i.e.
\begin{align*}
    \gamma^{(p)}_{\mu, \Lambda}(\mathcal{A}) &= \frac{1}{\sqrt{\text{det}(2\pi\Lambda)}}\int_{a \in \mathcal{A}}\exp\left\{-\frac{1}{2}({a}-{\mu})^\T\Lambda^{-1}({a}-{\mu})\right\}\;dx.
\end{align*}  Under Condition \ref{cond:asymlin}, $1-K_{T_n}\{T_n(Z),  \hat{P}\}$ is the usual asymptotically valid right-tail $p$-value for $\mathcal{H}_{\theta(P)}$ leveraging asymptotic normality.

\begin{proposition}\label{prop:gauss}
Suppose that $H_T(x, P)$ and $R_T(x, P)$ are continuous and strictly increasing and that Conditions \ref{cond:frac} - \ref{cond:sigma} hold. Then, Theorem \ref{thm:general} applies to the distributional estimator $K_{T_{n}}$ defined in (\ref{eq:asym}), such that the permutation distribution $R_{K_n}$ of $K_{T_n}\{T_n(Z),\hat{P}\}$ provides a test that is exact under $\mathcal{H}_P$ and asymptotically valid for $\mathcal{H}_{\theta(P)}$ under Conditions \ref{cond:frac}-\ref{cond:sigma}.   
\end{proposition}
The proof is presented in the appendix. Given Theorem \ref{thm:CLT}, it is a straightforward application of the continuous mapping theorem and Slutsky's theorem (both in the conventional case and for use within permutation distributions).
\begin{remark}\label{rem:MC}
For some choices of $T_n$ the transformation (\ref{eq:asym}) will have a familiar form, being the cumulative distribution function of a known distribution; examples are given in \S\S\ref{sec:bf}-\ref{sec:anova} and include normal and $\chi^2$ distributions. When this is not the case the probability required for the computation of (\ref{eq:asym}) may be replaced with a Monte-Carlo estimate from $B$ draws of a multivariate normal with mean zero and covariance $\hat{\Gamma}_\pi$. Importantly, using a Monte-Carlo estimate does not corrupt exactness of the test when $P_1=\ldots=P_k$.
\end{remark}

\subsection{Bootstrap prepivoting}\label{sec:boot}
The transformation outlined in \S \ref{sec:asym} makes explicit use of the limiting behavior for $T_{n}$, and requires the existence of a covariance estimator $\hat{\Sigma}$ satisfying Condition \ref{cond:sigma}. We now describe circumstances under which asymptotic validity for permutation tests when testing equality of parameters may be restored with the aid of bootstrap prepivoting \citep{ber87, ber88}. This avoids the need for an estimator $\hat{\Sigma}$, but will be shown to have appeal even when a convenient covariance estimator exists. This approach was first employed by \citet{chu16mult} to restore asymptotic validity for the permutation distribution of the max $t$-statistic under the null of equality of multivariate means. We now demonstrate that the approach is applicable to test statistics satisfying Condition \ref{cond:T} under additional regularity conditions on the involved estimators.

For a given permutation $\pi\in \cP_{n}$ and for $i=1,..,k$ let $Z^*_{\pi i1},\ldots,Z^*_{\pi in_i}$ be independent and identically distributed draws from $\hat{P}_{\pi i}$, the empirical distribution of $Z_{\pi i}$. Let $Z^*_{\pi i}$ be the $n_i \times {d}_P$ matrix whose $j$th row contains $Z^*_{\pi ij}$, and write $Z^*_{\pi}$ for the $n\times {d}_P$ matrix stacking $Z^*_{\pi 1}, ,\ldots,  Z^*_{\pi k}$ on top of one another. Write $\hat{\Theta}^*_\pi = \hat{\Theta}(Z_\pi^*)$,  $\hat{C}^*_\pi=\hat{C}(Z_\pi^*)$, and $\hat{\eta}^*_\pi = \hat{\eta}(Z_\pi^*)$. 

Let $\breve{T}_{n}$ be the bootstrap modification of $T_n$, defined as
\begin{align}\label{eq:Tboot}
\breve{T}_n(Z^*_\pi, Z_\pi) &= g\{n^{1/2}(\hat{\Theta}^*_\pi - \hat{\Theta}_\pi)\hat{C}^*_\pi, \hat{\eta}^*_\pi\},
\end{align}and define the bootstrap distribution function $J_{T_{n}}(x, \hat{P}_\pi)$ as
\begin{align}\label{eq:bootdist}
J_{T_n}(x, \hat{P}_\pi) = \text{pr}\{\breve{T}_n(Z_\pi^*, Z_\pi) \leq x\mid Z_\pi\}.
\end{align}
The modification $\breve{T}_n$ is required for centering the bootstrap null distribution. The bootstrap prepivoted test statistic is then $J_{T_n}\{T_n(Z_\pi), \hat{P}_\pi\}$. For the bootstrap to be applicable in our setup, additional conditions are required.
\begin{condition}\label{cond:asymlinboot}
Let $X^*_{i1},\ldots,X^*_{in_i}$ be independent and identically distributed draws from $\hat{P}_i$, and let $X^*_i$ be the $n_i\times {d}_P$ matrix whose $j$th row is $X^*_{ij}$. Then, for $i=1,\ldots,k$,
\begin{align*}
n_i^{1/2}\{\hat{\theta}_\ell(X^*_i) - \hat{\theta}_\ell(X_i)\} &= n_i^{-1/2}\sum^{n_i}_{j=1}\{{f}_{P_{i}, \ell}(X^*_{ij}) - {f}_{P_{i}, \ell}(X_{ij})\} + \mathcal{E}_{in_i},
\end{align*} 
where $\mathcal{E}_{in_i}$ satisfies $\P(n_i^{1/2} \mathcal{E}_{in_i} \geq \epsilon \mid X_i)\rightarrow 0$ in probability as $n_i\rightarrow \infty$ for any $\epsilon > 0$. Furthermore, let $Y_1,\ldots,Y_n$ be distributed according to the mixture $\bar{P} = \sum_{i=1}^kP_i$, and let $Y^*_1,\ldots,Y^*_n$ be independent and identically distributed draws from the empirical distribution of $Y_1,\ldots,Y_n$, with $Y$ and $Y^*$ being $n\times {d}_P$ matrices whose $j$th rows are $Y_j$ and $Y^*_j$.  Then, 
\begin{align*}
n^{1/2}\{\hat{\theta}_\ell(Y^*) - \hat{\theta}_\ell(Y)\} &= n^{-1/2}\sum^{n}_{j=1}\{{f}_{\bar{P}, \ell}(Y^*_{j}) - {f}_{\bar{P}, \ell}(Y_{j})\} + \bar{\mathcal{E}}_{n},
\end{align*} 
where $\bar{\mathcal{E}}_{n}$ satisfies $\P(n^{1/2} \bar{\mathcal{E}}_{n} \geq \epsilon \mid Y)\rightarrow 0$ in probability as $n\rightarrow \infty$ for any $\epsilon > 0$. 
\end{condition}
\begin{condition}\label{cond:Tboot} $\hat{\eta}^*$ and $\hat{\eta}^*_\Pi$ satisfy
\begin{align*}
\P\{|\hat{\eta}(Z^*) - \eta| > \epsilon \mid Z\} \rightarrow 0 ;\;\; \P\{|\hat{\eta}(Z^*_\Pi) - \bar{\eta}| > \epsilon \mid Z_\Pi\} \rightarrow 0
\end{align*} in probability as $n\rightarrow \infty$ for any $\epsilon > 0$. The analogous statements hold for $\hat{C}(Z^*)$ and $\hat{C}(Z^*_\Pi)$ with limiting values $C$ and $\bar{C}$ respectively. \end{condition}
\begin{proposition}

\label{prop:boot} Suppose that $H_T(x, P)$ and $R_T(x, P)$ are continuous and strictly increasing and that Conditions \ref{cond:frac}-\ref{cond:T} and \ref{cond:asymlinboot}-\ref{cond:Tboot} hold. Then, Theorem \ref{thm:general} applies to  $J_{T_{n}}$ defined in (\ref{eq:bootdist}).
\end{proposition} 
The proof of this proposition closely mirrors that of Theorem 2.6 of \citet{chu16mult}. The result for the permutation distribution $R_{J_n}$ is proven in the appendix, while the result for the true unconditional distribution $H_{J_n}$ is standard under these conditions and omitted. \citet{liu89} describe the representation given in Condition \ref{cond:asymlinboot} and sufficient conditions giving rise to it. Of particular interest, Hadamard differentiability of $\theta(\cdot)$ at $P_1,\ldots,P_k$ and $\bar{P}$ is sufficient for Condition \ref{cond:asymlinboot} to hold. Verifying Condition \ref{cond:Tboot} for mean-like statistics often involves applying certain triangular array weak laws, such as that of \citet[Theorem 2.2, part b]{bic81} or \citet[Lemma 15.4.1]{leh05}. Lemma 4.1 of \citet{dum13} shows that Condition \ref{cond:Tboot} may equivalently be expressed in terms of the unconditional convergence in probability of $\hat{\eta}(Z^*)$, $\hat{\eta}(Z_\Pi^*)$, $\hat{C}(Z^*)$, and $\hat{C}(Z^*_\Pi)$, accounting for randomness in both the observed and the bootstrap sample for $Z^*$ and observed sample, the bootstrap sample, and the permutation draw for $Z^*_\Pi$; see also \citet{buc19}.

In practice both the bootstrap law $J_{T_n}$ and the permutation distribution $R_{J_n}$ will be approximated through Monte-Carlo simulation. In \S B of the appendix we provide pseudocode illustrating the general implementation. We also provide \texttt{R} code for the test statistics described in \S\ref{sec:sim} and in \S A of the appendix.

\section{Higher-order corrections}\label{sec:higher}
\subsection{Applying bootstrap prepivoting to pivotal and non-pivotal quantities}
Given the additional computational burden imposed by embedding a bootstrap scheme within a permutation test, the approach of \S \ref{sec:boot} may appear unappealing relative to that of \S \ref{sec:asym}.  We now illustrate that higher order refinements may be attained for prepivoted permutation tests under the equality of parameters through bootstrap prepivoting, so long as the statistic being bootstrapped admits an asymptotic expansion whose first term is pivotal. We then describe how Gaussian prepivoting may be used in concert with bootstrap prepivoting to furnish an expansion with a pivotal leading term.

Mirroring \S 3 of \citet{ber88}, in what follows we assume the existence of certain asymptotic expansions for the unconditional, permutation, and bootstrap distributions of test statistics $T_{n}$ satisfying Condition \ref{cond:T}. Rather than providing primitive conditions under which the expansions exist for specific test statistics such as those of the form specified in Condition \ref{cond:T}, in our exposition we will highlight interesting consequences of the expansions assuming their existence. Generally the expansions hold under the smooth function model of \citet{bha78}; see also \citet[\S 2.4]{hal92}. Various forms of differentiability such as Hadamard differentiability which imply bootstrap consistency are generally insufficient for the existence of these expansions. In \S \ref{sec:bf}, we provide sufficient conditions under which the expansions exist in the case of the univariate difference in means.

\begin{condition}\label{cond:expand} $H_{T_{n}}$ and $R_{T_{n}}$ admit asymptotic expansions of the following form for some integer $j\geq 1$:
\begin{align}
\label{eq:H}H_{T_{n}}(x, P) &= A_{T_n}(x, P) + n^{-j/2}a_{T_n}(x, P) + O(n^{-(j+1)/2});\\
\label{eq:R}R_{T_{n}}(x,\hat{P}) &= B_{T_n}(x,\hat{P}) + n^{-j/2}b_{T_n}(x,\hat{P}) + O_P(n^{-(j+1)/2}),
\end{align}uniformly in $x$, where the functions on the right-hand side are continuous in $x$. 

\end{condition}
\begin{condition}\label{cond:dhat}
${K}_{T_n}$ in (\ref{eq:asym}) satisfies
\begin{align*}
{K}_{T_n}(x,\hat{P}) - A_{T_n}(x, P) = O_P(n^{-1/2});&\;\;\; {K}_{T_n}(x,\hat{P}_\Pi) - B_{T_n}(x,  \hat{P}) = O_P(n^{-1/2})\end{align*}
uniformly in $x$. Moreover, uniformly in $x$, there exists an estimator ${k}_{T_n}(x,\hat{P})$ satisfying
\begin{align*}
{k}_{T_n}(x,\hat{P}) - a_{T_n}(x, P) = O_P(n^{-1/2});&\;\;\; {k}_{T_n}(x,\hat{P}_\Pi) - b_{T_n}(x,  \hat{P}) = O_P(n^{-1/2}).
\end{align*} Both $K_{T_n}$ and $k_{T_n}$ are continuous in $x$.\end{condition}
\begin{condition}\label{cond:expandboot}
The bootstrap estimator $J_{T_n}$ for the distribution of $T_{n}$ admits the following expansions for some integer $j\geq 1$ uniformly in $x$:
\begin{align}
\label{eq:Hboot} J_{T_{n}}(x,\hat{P}) &= K_{T_n}(x,\hat{P}) + n^{-j/2}k_{T_n}(x,\hat{P}) + O_P(n^{-(j+1)/2});\\
\label{eq:Rboot} J_{T_{n}}(x,\hat{P}_\Pi) &= K_{T_n}(x,\hat{P}_\Pi) + n^{-j/2}{k}_{T_n}(x,\hat{P}_\Pi) + O_P(n^{-(j+1)/2}).
\end{align} \end{condition}



Let $H_{J_n}$ and $R_{J_n}$ be the unconditional and permutation distributions for the bootstrap prepivoted $J_{T_n}\{{T}_n(Z), \hat{P}\}$. The following theorems provide the order of errors in rejection probability using a bootstrap prepivoted permutation test under $\mathcal{H}_{\theta(P)}$. 
\begin{theorem}\label{thm:dep}
Suppose that Conditions \ref{cond:T} and \ref{cond:expand} - \ref{cond:expandboot} hold. Consider replacing the test statistic $T_{n}(Z_\pi)$ with the bootstrap prepivoted $J_{T_{n}}\{T_{n}(Z_\pi),  \hat{P}_\pi\}$. Let $U(x)$ be the distribution function of a uniform random variable on $[0,1]$. Then, uniformly over $x$, 
\begin{align*}
H_{J_{n}}(x, P)= U(x) + O(n^{-1/2});\;\; R_{J_{n}}(x,\hat{P}) = U(x) + O_P(n^{-1/2}),
\end{align*}
 As a result, $H_{J_{n}}(x,P) - R_{J_{n}}(x,\hat{P}) = O_P(n^{-1/2})$ uniformly in $x$.
\end{theorem}
\begin{theorem}\label{thm:indep}
Consider the setup of Theorem \ref{thm:dep}, but strengthen Condition \ref{cond:dhat} so that $K_{T_n}(x, \hat{P}_\Pi) = K_{T_n}(x)$, $A_{T_n}(x,P) = A_{T_n}(x),$ and  $B_{T_n}(x,  \hat{P}) = B_{T_n}(x)$, with $A_{T_n}(x) = B_{T_n}(x) = K_{T_n}(x)$. Then,
\begin{align*}
H_{J_{n}}(x, P)= U(x) + O(n^{-(j+1)/2});\;\; R_{J_{n}}(x,\hat{P}) = U(x) + O_P(n^{-(j+1)/2}),
\end{align*}
such that $H_{J_{n}}(x, P) - R_{J_{n}}(x,\hat{P}) = O_P(n^{-(j+1)/2})$ uniformly in $x$.
\end{theorem}The proofs of Theorems \ref{thm:dep} and \ref{thm:indep} are presented in the appendix, and closely follow the presentation of Cases 1 and 2 in  \citet[\S 3.1]{ber88}. In \S C of the appendix, we discuss how repeated applications of bootstrap prepivoting may be used to further improve the order of the error in rejection probability. If the asymptotic expansions exist with the required number of terms, $\ell\geq 1$ applications of bootstrap prepivoting would yield an error order of $O(n^{-(j_1 + \ell-1)/2})$, where $j_1 = 1$ under the conditions of Theorem \ref{thm:dep} and $j_1=j+1$ under the conditions of Theorem \ref{thm:indep}. The required layers of bootstrapping become computationally burdensome as the number of iterations increase while providing diminishing improvements to the order of approximation, making anything beyond two applications of bootstrap prepivoting within a permutation test unappealing in most settings.   

\subsection{Bootstrap-after-Gaussian prepivoting}\label{sec:bg}

Assuming that the expansions are valid, Theorems \ref{thm:dep} and \ref{thm:indep} show that whether or not the bootstrap furnishes a higher-order correction depends upon whether or not $A_{T_n} = B_{T_n} = K_{T_n}$. Under the setup of these theorems the Gaussian prepivoted test statistic $K_{T_n}\{T_n(Z), \hat{P}\}$ results in a test with an error in rejection probability $O(n^{-1/2})$, the same order recovered through Theorem \ref{thm:dep}. As a result, bootstrap prepivoting can only provide an improvement in the order of the error in rejection probability over Gaussian prepivoting when $T_n$ is an asymptotic pivot.

We now show that the approaches of \S\ref{sec:asym} and \S\ref{sec:boot} may be combined to ensure the equality of the lead terms in the expansions, placing us within the favorable scenario of Theorem \ref{thm:indep}. Suppose in addition to Condition \ref{cond:sigma} that the estimator $\hat{\Gamma}(\cdot)$ used in the Gaussian prepivoting transformation satisfies Condition \ref{cond:Tboot} for $\hat{\Gamma}(Z^*)$ and $\hat{\Gamma}(Z^*_\Pi)$. Consider using bootstrap prepivoting on a Gaussian prepivoted test statistic, yielding the bootstrap distribution function
\begin{align*}
J_{K_n}(x, \hat{P}_\pi) &= \text{pr}\left[K_{T_n}\{\breve{T}_n(Z_\pi^*, Z_\pi), \hat{P}_\pi^*\} \leq x\mid Z_\pi\right],
\end{align*} 
where $\hat{P}^*_\pi = \prod_{i=1}^k\prod_{j=1}^{n_i}\hat{P}^*_{\pi i}$ and $\hat{P}^*_{\pi i}$ is the empirical distribution of $Z^*_{\pi i}$. The corresponding prepivoted test statistic is $J_{K_n}[K_n\{T_n(Z_\pi),\hat{P}_\pi\}, \hat{P}_\pi]$. After Gaussian prepivoting, by Theorem \ref{prop:gauss} the lead term of the expansions for both the unconditional and permutation distributions of $K_{T_n}\{T_n(Z),  \hat{P}\}$ are $A_{K_n}(x, P) = B_{K_n}(x,  \hat{P}) = U(x)$, the distribution function for the standard uniform. Meanwhile,  $K_{T_n}\{T_n(Z_\pi),\hat{P}_\pi\}$ as defined in (\ref{eq:asym}) is itself a test statistic of the form (\ref{eq:T}) whenever $T_n(Z_\pi)$ is, and $K_{T_n}\{\breve{T}_n(Z_\pi^*, Z_\pi),\hat{P}^*_\pi\}$ is of the form (\ref{eq:Tboot}) whenever $\breve{T}_n(Z_\pi^*, Z_\pi)$ is; see the proof of Proposition \ref{prop:gauss} for details. Therefore, by Proposition \ref{prop:boot} the lead term in the relevant bootstrap expansions would also be $U(x)$ in regular cases. So long as the expansions exist Theorem \ref{thm:indep} applies when bootstrapping a Gaussian prepivoted statistic.

As previously discussed Gaussian prepivoting requires a suitable covariance estimator $\hat{\Sigma}(\cdot)$, used to form $\hat{\Gamma}(Z_\pi)$ in (\ref{eq:asym}). If a natural covariance estimator is unavailable, applying two layers of bootstrap prepivoting can also provide a higher-order correction; see \S C of the appendix for details.

\subsection{Analytical alternatives}  

In situations where expressions for the expansions (\ref{eq:H}) - (\ref{eq:R}) are known in closed form, permutation tests with higher-order correctness under $\mathcal{H}_{\theta(P)}$ can be attained using suitable estimates for the terms in the expansions. This can avoid the computational costs of embedding the bootstrap within a permutation test, or can be used in conjunction with the bootstrap to further improve the order of the error in rejection probability. The approach closely follows the developments in \S 3.2 of \citet{ber88}. 

Consider the modified test statistic
\begin{align}\label{eq:edge}E_{T_n}\{T_{n}(Z_\pi),\hat{P}_\pi\} &=  {K}_{T_n}\{T_{n}(Z_\pi),\hat{P}_\pi\} + n^{-j/2}{k}_{T_n}\{T_{n}(Z_\pi),\hat{P}_\pi\},\end{align}
with $K_{T_n}$ as in (\ref{eq:asym}) and $k_{T_n}$ an estimator satisfying Condition \ref{cond:dhat}. Let $H_{E_n}$ and $R_{E_n}$ be the unconditional and permutation distributions of $E_{T_n}\{T_n(Z_\pi), \hat{P}_\pi\}$. The following theorems again show that the order of the difference between  $H_{E_n}(x, P)$ and $R_{E_n}(x,  \hat{P})$ depends upon whether $K_{T_n}(x, \cdot)$ aligns with the lead terms in (\ref{eq:H}) and (\ref{eq:R}).

\begin{theorem}\label{thm:dep2}
Suppose that Conditions \ref{cond:T}, \ref{cond:expand} and \ref{cond:dhat} hold. Then, uniformly in $x$,
\begin{align*}
H_{E_{n}}(x,P) - R_{E_{n}}(x,\hat{P}) = O_P(n^{-1/2}).
\end{align*}
\end{theorem}

\begin{theorem}\label{thm:indep2}
Suppose that Conditions \ref{cond:T}, \ref{cond:expand} and \ref{cond:dhat} hold. Suppose further that $K_{T_n}(x, \hat{P}_\Pi) = K_{T_n}(x)$, $A_{T_n}(x, P) = A_{T_n}(x)$ in (\ref{eq:H}), and $B_{T_n}(x,  \hat{P}) = B_{T_n}(x)$ in (\ref{eq:R}), with $K_{T_n}(x) = A_{T_n}(x) = B_{T_n}(x)$. Then, uniformly in $x$,
\begin{align*}
H_{E_n}(x, P) - R_{E_n}(x,  \hat{P}) &= O_P(n^{-(j+1)/2}).\end{align*} 
\end{theorem}The proofs of Theorems \ref{thm:dep2} and \ref{thm:indep2} are analogous to those of Theorem \ref{thm:dep} and \ref{thm:indep} respectively and are omitted; simply replace $J_{T_n}$ with $E_{T_n}$, $H_{J_n}$ with $H_{E_n}$, and $R_{J_n}$ with $R_{E_n}$ in the proofs. Theorem \ref{thm:indep2} shows the potential benefit of adding an estimate of the next term of the expansion of $T_n$ to a Gaussian prepivoted test statistic whenever the base statistic $T_n$ is already pivotal: an improvement in the order of the error rejection probability may be attained without using the bootstrap. Moreover, one may apply bootstrap prepivoting to $E_{T_n}$ in order to attain the same higher-order accuracy as applying a double bootstrap (so long as the required expansions exist for $E_{T_n}$). Theorem \ref{thm:dep2} illustrates that when $T_n$ is not pivotal, there is no benefit from adding an estimate of the subsequent term for the order of the error in rejection probability relative to simply using the Gaussian prepivoted test statistic, much as there was no benefit to applying bootstrap prepivoting in this case. Under Conditions \ref{cond:expand} and \ref{cond:dhat} and assuming the setup of Theorem \ref{thm:dep2} the Gaussian prepivoted test statistic ${K}_{T_n}\{T_n(Z_\pi), \hat{P}_\pi\}$ also satisfies $H_{K_n}(x, P) - R_{K_n}(x,\hat{P}) = O_P(n^{-1/2})$.

\section{Examples and illustrations}\label{sec:sim}
\subsection{Notation for the examples}\label{sec:notation}
We present simulation studies to illustrate the benefits of prepivoted permutation tests. We will focus on properties of these tests under $\mathcal{H}_{\theta(P)}$, highlighting both asymptotic correctness and higher-order accuracy for various prepivoted permutation tests while illustrating the impropriety of other seemingly natural permutation tests. The primary benefit of using prepivoted permutation tests over alternative large-sample tests for $\mathcal{H}_{\theta(P)}$ such as bootstrap hypothesis tests is the exactness of inference for any $n$ under equality in distribution. While we do not showcase this property in simulation studies, it in combination with the theoretical results and simulation studies presented here provide compelling motivation for deploying prepivoted permutation tests in practice.

Sections \ref{sec:bf} - \ref{sec:anova} concern inference for equality of population means, denoted $\mathcal{H}_{\mu(P)}$, with scenarios varying in whether the two- or $k$-sample problem is under consideration and whether the response is univariate or multivariate. In these examples, we will let $\hat{\mu}_{\pi i} = \hat{\mu}(Z_{\pi i})$ be the sample mean of $Z_{\pi i}$, with $\hat{\mu}_{\pi} = (\hat{\mu}_{\pi 1},..,\hat{\mu}_{\pi k})^\T$; $\hat{\mu}_{\pi i}$ may be a scalar or vector-valued depending upon the context. $\hat{\Sigma}_{\pi i}$ will be either the sample variance or the sample covariance matrix for $Z_{\pi i}$, and $\hat{\Gamma}_\pi$ will be the $kd\times kd$ block-diagonal matrix with $(n/n_i)\hat{\Sigma}_{\pi i}$ in the $i$th of $k$ blocks. Test statistics will involve one of two contrast matrices. First, let $C_1$ be the $k\times k(k-1)/2$ matrix of contrasts yielding all pairwise differences between the $k$ group means. For instance, when $k=2$ $C_1 = (1,-1)^\T$; for $k=3$, $C_1$ is $3\times3$  with columns $(1,-1,0)^\T, (1,0,-1)^\T$, and $(0, 1, -1)^\T$. Second, let $b = (n_1/n,\ldots,n_k/n)^\T$, let $B$ be the $k\times k$ matrix with $B_{ij} = b_i$ for all $i$, and define $C_2 = I_{k\times k} - B$. Finally, a few common hypothesis tests for $\mathcal{H}_{\mu(P)}$ assume homogeneity of variances, with the corresponding statistics using the usual pooled estimator of the variance. Through our examples we will show that bootstrap and Gaussian prepivoting may be applied to these statistics to restore asymptotic validity even when variances are heterogeneous. To define the pooled variance let $A$ be the $n\times k$ matrix whose $i$th column indicates membership in the $i$th group, such that $A_{ji} = 1$ if $j\in \mathcal{I}_i$, 0 otherwise. Let $H=A(A^\T A)^{-1}A^\T$,  define $\hat{\nu}_\pi = Z_\pi^\T(I_{n\times n} - H)Z_\pi/(n-k)$ as the pooled variance estimator, and let $\hat{\Lambda}_\pi$ be the $kd\times kd$ block-diagonal matrix with $(n/n_i)\hat{\nu}_\pi$ in the $i$th of $k$ blocks. In these examples the required convergence in probability of $\hat{\nu}(Z_\Pi)$, $\hat{\Lambda}(Z_\Pi)$ and $\hat{\Gamma}(Z_\Pi)$ under Condition \ref{cond:T} may be verified using Lemma 5.3 of \citet{chu13}.

\subsection{Nonparametric Behrens-Fisher}\label{sec:bf}
We begin with the familiar two-sample, univariate case, and compare the performance of various permutation tests when only equality of expectations holds. The candidate test statistics are
\begin{align*}T_{n}(Z_\pi) &= n^{1/2}\hat{\mu}_\pi^\T C_1 = n^{1/2}(\hat{\mu}_{\pi 1} - \hat{\mu}_{\pi 2});\\
S_{n}(Z_\pi) &= \frac{T_{n}(Z_\pi)}{\{C_1^T\hat{\Gamma}_\pi C_1\}^{1/2}} = \frac{n^{1/2}(\hat{\mu}_{\pi 1} - \hat{\mu}_{\pi 2})}{\{(n/n_i)\hat{\Sigma}_{\pi 1} + (n/n_2)\hat{\Sigma}_{\pi 2}\}^{1/2}};\\
K_{T_{n}}\{T_n(Z_\pi),\hat{P}_\pi\} &= K_{S_n}\{S_n(Z_\pi), \hat{P}_\pi\} = \Phi\left\{S_n(Z_\pi)\right\};\\
{E}_{S_n}\{S_n(Z_\pi),\hat{P}_\pi\} &= \Phi\left\{S_n(Z_\pi)\right\} + \frac{1}{6} \left(\frac{\hat{\xi}(Z_\pi)}{\left\{C_1^T\hat{\Gamma}_\pi C_1\right\}^{3/2}}\right)\phi\{S_n(Z_\pi)\}\left[2\{S_n(Z_\pi)\}^2 + 1\right]n^{-1/2},
\end{align*}
where
$\hat{\xi}(Z_\pi) = (n/n_1)^2\hat{\kappa}_{\pi 1} - (n/n_2)^2\hat{\kappa}_{\pi 2}$, $\hat{\kappa}_{\pi i}$ is the  sample central third moment for group $i$,  and  $\Phi(\cdot)$ and $\phi(\cdot)$ are the cumulative distribution function and density function for the standard normal distribution. 

$T_n$ is the difference in means. From Theorem \ref{thm:CLT},  this test statistic may yield anti-conservative inference under $\mathcal{H}_{\mu(P)}$. $K_{T_n}$ applies the approach of \S\ref{sec:asym} and is the simplification of (\ref{eq:asym}) for this choice of $T_n$. $S_n$ is the usual $t$-statistic studentized using the unpooled variance estimator, advocated in \citet{chu13} to restore first-order accuracy to permutation tests under $\mathcal{H}_{\mu(P)}$. Observe that $K_{T_n}$ is a monotone increasing function of $S_n$, such that permutation tests using $S_n$ and $K_{T_n}$ will yield identical $p$-values. Due to this equivalence the success of studentization for restoring first-order accuracy can be understood in light of Proposition \ref{prop:gauss}.  $E_{S_n}$ applies (\ref{eq:edge}) to $S_n$, with the second term in $E_{S_n}$ representing an estimate of the second term in the Edgeworth expansion for $S_n$ \citep[\S 4]{abr85}. The estimate satisfies Condition \ref{cond:dhat}, and the performance of $E_{S_n}$ will reflect Theorem \ref{thm:indep2} due to pivotality of the lead term in the expansions for $S_n$. 

Sufficient conditions for Conditions \ref{cond:expand} and \ref{cond:expandboot} to hold for $T_n$ are that $E||X_{i1}||^{4+\delta} < \infty$ for $i=1,2$ for some $\delta > 0$ and that $P_i$ is continuous for at least one $i$. For $S_n$ and $K_{T_n}$, $E||X_{i1}||^{8+\delta} < \infty$ along with the aforementioned continuity assumption are sufficient, while $E||X_{i1}||^{10+\delta} < \infty$ and the continuity assumption are sufficient for the expansions for $E_{S_n}$ to satisfy Conditions \ref{cond:expand} and \ref{cond:expandboot}. Validity of the expansions (\ref{eq:H}) and (\ref{eq:Hboot}) follows from classical arguments in the $k$-sample case under these assumptions; see \citet{bab83}, \citet{hal88two}, and \citet{hal92} among many. The expansions for the randomization distributions (\ref{eq:R}) and (\ref{eq:Rboot}) are less conventional but can nonetheless be established using existing results. Validity of (\ref{eq:Rboot}) can be established through various approaches, such as through results from \citet{liu89} for the bootstrap under non independent and identically distributed models. Interestingly, we can alternatively use the contiguity results in \citet[Lemma 5.3]{chu13} to leverage results on Edgeworth expansions for the bootstrap in the $k$-sample problem in the special case that all $k$ samples arise as independent and identically distributed draws from a common distribution $\bar{P}$. Essential for applying these contiguity results is the fact that the moment conditions and continuity requirements on the individual populations are all that is required for the conventional $k$-sample case. For (\ref{eq:R}), one can use results from finite population survey sampling under the assumption that the finite population is itself generated from a superpopulation model. For instance, both $T_n$ and $S_n$ may be expressed in the form required for Theorem 2 of \citet{bab85}; see also \citet{boo93}, \citet{boo94} and \citet{bab96} for related results. In particular, results in \citet{bab96} provide for the existence of the expansions with remainder $O_P(n^{-j/2})$ for $j>1$ under our setup.     

For each test we will reject for large values of the test statistic, with desired Type I error rate of $\alpha = 0.05$. For each test statistic we further consider the permutation distributions after bootstrap prepivoting, yielding 8 test statistics in total. In each simulation scenario $X_{11},\ldots,X_{1n_1}$ will be independent and identically distributed as $-(\mathcal{E}(1/8) - 8)$ and $X_{21},\ldots,X_{2n_2}$ will be independent and identically distributed as $\mathcal{E}(1/5)-5$, where $\mathcal{E}(\lambda)$ is the exponential distribution with rate $\lambda$. The sample sizes $n_1$ and $n_2$ will be varied across scenarios, with $n_1\neq n_2$. In each scenario we simulate 10,000 data sets. Permutation tests within each data set are based upon 999 draws from the permutation distribution. For the bootstrap prepivoted test statistics, 500 bootstrap samples are drawn for each permutation.

The first row of Table \ref{tab:sim1} shows the expected accuracy of the permutation distributions. Permuting $T_n$ itself will be asymptotically invalid because it is not asymptotically pivotal. Theorem \ref{thm:dep} applies to the bootstrap prepivoted transform of $T_n$, while Theorem \ref{thm:indep} applies to the bootstrap prepivoted transforms of the remaining test statistics. For $S_n$ and $K_{S_n}$, the order of the error in rejection probability is $O(n^{-1/2})$ before applying bootstrap prepivoting, while the rate is $O(n^{-1})$ for $E_{S_n}$ before applying bootstrap prepivoting by Theorem \ref{thm:indep2}.

\begin{table}
\caption{Order of the errors in rejection probability (ERP) and simulated Type I error rates for permutation tests under $\mathcal{H}_{\mu(P)}$ in the univariate, two-sample setting. For each simulation $n_1 = 0.3 n$, and the nominal level of the tests is $\alpha = 0.05$. $T_n$ and $S_n$ are the unstudentized and studentized difference in means respectively. $K_{T_n}$ applies (\ref{eq:asym}) to $T_n$, and $E_{S_n}$ applies (\ref{eq:edge}) to $S_n$. For any given sample size $n$, the maximal estimated standard error for the difference between two estimated rejection probabilities is 0.0024.}\label{tab:sim1}
\begin{center}
\begin{tabular}{l c c c c c c c c}
&\multicolumn{4}{c}{Original} &\multicolumn{4}{c}{Bootstrap Prepivoted}\\
&$T_n$&$S_n$&$K_{T_n}$&$E_{S_n}$& $T_n$&$S_n$&$K_{T_n}$&$E_{S_n}$\\
Error Order & $O(1)$ & $O(n^{-1/2})$ & $O(n^{-1/2})$ & $O(n^{-1})$& $O(n^{-1/2})$&$O(n^{-1})$&$O(n^{-1})$&$O(n^{-3/2})$\\
$n=50$ & 0.125 & 0.119 & 0.119 & 0.087 & 0.130 & 0.075&0.075 & 0.068\\
$n=100$ & 0.105 & 0.090 & 0.090 & 0.061 & 0.098 & 0.055&0.055 & 0.052\\
$n=250$ &0.105 & 0.079 & 0.079 & 0.055 & 0.083 & 0.050 &0.050 & 0.049\\
$n=1000$&0.092 & 0.059 & 0.059 & 0.049 & 0.061 & 0.047 & 0.047 & 0.048\\
$n=5000$&0.092 & 0.054 & 0.054 & 0.050 & 0.055 & 0.048 & 0.048& 0.051\\
\end{tabular}
\end{center}
\end{table}


Table \ref{tab:sim1} shows the results of our simulation study with various choices of $n_1$ and $n_2$. Even with the largest values of $n_1$ and $n_2$, $T_n$ does not provide a rejection rate close to the nominal level. This can be understood in light of Theorem \ref{thm:CLT} through a comparison of the covariances (\ref{eq:cov}) and (\ref{eq:covperm}): the unconditional limit distribution of $T_n$ has a larger variance than that of the probability limit of $R_{T_n}$, inflating the rejection rate even asymptotically. For the remaining 7 test statistics, for $n_1$ and $n_2$ sufficiently large the rejection rates approach the nominal level, a reflection of Propositions \ref{prop:gauss} and \ref{prop:boot}. For smaller values of $n$,  for each test statistic bootstrap prepivoting brings the rejection rate closer to the nominal level. Furthermore,  both before and after applying bootstrap prepivoting $E_{S_n}$ comes closer to the nominal error rate than $K_{T_n}$ or $S_n$, which are in turn closer to the nominal rate than $T_n$. This reflects the developments in \S \ref{sec:higher}. As desired, these higher-order improvements are particularly noticeable at small values of $n$. At $n=50$  the bootstrap prepivoted version of $E_{S_n}$ has an estimated Type I error rate of 0.068 compared to 0.119 for $S_n$ without bootstrap prepivoting. At $n=100$, these values are 0.052 and 0.090 respectively. Even at $n=1000$, the permutation test based upon $S_n$ has a rejection rate slightly above the nominal level despite asymptotic validity.


The permutation tests based using $E_{S_n}$ and its bootstrap prepivoted transformation provide appealing approaches to the nonparametric Behrens-Fisher problem. The permutation test using $E_{S_n}$ is exact for any sample size if $P_1=P_2$, and provides an error in rejection probability of order $O(n^{-1})$ under $H_{\mu(P)}$. After applying bootstrap prepivoting to $E_{S_n}$, the order in the error in rejection probability for the resulting permutation test is reduced to $O(n^{-3/2})$ while maintaining exactness if $P_1=P_2$. Through additional applications of bootstrap prepivoting to $E_{S_n}$ the error in rejection probability can be further reduced under $H_{\mu(P)}$ while maintaining exactness when $P_1=P_2$; see the discussion of bootstrap iteration in \S C of the appendix for more details.
\subsection{Multivariate tests: Hotelling and max absolute $t$}\label{sec:mult}
In this section we keep $k=2$ but now consider multivariate responses. Three common test statistics for testing $\mathcal{H}_{\mu(P)}$ in this situation are the pooled and unpooled Hotelling statistics along with the max absolute $t$ statistic ($L_n$, $S_n$, and $M_n$ respectively):
\begin{align*}
L_n(Z_\pi) &= n(\hat{\mu}_\pi C_1)^\T(C_1^\T\hat{\Lambda}_\pi C_1)^{-1}(\hat{\mu}_\pi C_1);\\
S_n(Z_\pi) &= n(\hat{\mu}_\pi C_1)^\T(C_1^\T\hat{\Gamma}_\pi C_1)^{-1}(\hat{\mu}_\pi C_1);\\
M_n(Z_\pi) &= \underset{1\leq \ell \leq d}{\max}\;\;\frac{n^{1/2}|(\hat{\mu}_\pi C_1)_\ell|}{\{(C_1^\T\hat{\Gamma}_\pi C_1)_{\ell \ell}\}^{1/2}}.
\end{align*} 
In addition to investigating the performance of these three test statistics under $\mathcal{H}_{\mu(P)}$, we consider the impact of (a) Bootstrap prepivoting the original statistics; (b) Gaussian prepivoting the original statistics; and (c) Bootstrap prepivoting the Gaussian prepivoted versions of the statistics, yielding 12 total test statistics. \citet{chu16mult} considered the use of the unpooled version of Hotelling's statistic $S_n(Z_\pi)$ along with the bootstrap prepivoted version of the max absolute $t$ statistic $J_{M_n}\{M_n(Z_\pi), \hat{P}_\pi\}$ for inference on $\mathcal{H}_{\mu(P)}$ using permutation tests. They further demonstrated that without further modification, $L_n(Z_\pi)$ cannot be used for inference for $\mathcal{H}_{\theta(P)}$.

Applying Gaussian prepivoting to the unpooled Hotelling statistic yields the test statistic $G_{d}\{S_n(Z_\pi)\}$, where $G_{df}$ is the distribution function of a $\chi^2$ random variable with $df$ degrees of freedom. As this is a monotone increasing transformation of $S_n$, permutation $p$-values attained after Gaussian prepivoting will be identical to those attained without Gaussian prepivoting. Through this equivalence, Proposition \ref{prop:gauss} proves the asymptotic validity of permutation inference using $S_n$ for testing equality in means. Applying Gaussian prepivoting to $L_n(Z_\pi)$ requires computing the distribution function for general quadratic forms of multivariate Gaussians; we use the algorithm of \citet{dav80}, implemented by the function \texttt{davies} within the \texttt{R} package \texttt{CompQuadForm} to accomplish this. Applying Gaussian prepivoting to $M_n(Z_\pi)$ requires the computation of the probability that correlated mean-zero Gaussians fall within a rectangular region symmetric about the origin. While this can be accomplished using the \texttt{pmvnorm} function in \texttt{R}, for improved speed we proceed along the lines of Remark \ref{rem:MC}, drawing $B=1000$ Monte-Carlo draws from a multivariate normal with mean zero and covariance $\hat{\Gamma}_\pi$ in order to approximate (\ref{eq:asym}) within each permutation. Bootstrap prepivoting simply computes the bootstrap distribution of (\ref{eq:Tboot}), with the function $g$ within Condition \ref{cond:T} varying with the test statistic being employed. 

Our response variables are $d_P=15$ dimensional. For  $j=1,\ldots,n_1$ we draw $X_{1j}$ independently and identically distributed as $\exp\{\mathcal{N}_{15}(0_{15}, I_{15\times 15})\}$, where $\mathcal{N}_{15}$ is the multivariate normal distribution and $0_{15}$ is a vector of length 15 with zeroes in each element.  For $j=1,\ldots,n_2$ we then draw $X_{2j}$ independently and identically distributed as $\mathcal{N}_{15}(\mu, V)$, where $V_{\ell \ell} = 1$, $V_{\ell \ell'} = 0.8$ for $\ell\neq \ell'$, and $\mu_{\ell} = \exp(0.5)$ for $\ell,\ell' = 1,\ldots,15$. We set $n_1 = 0.3n$ and $n_2=0.7n$. We simulate 5000 data sets using the above generative model with different values for $n$. For each data set, permutation inference proceeds using 999 permutations. When bootstrap prepivoting is used, we draw 200 bootstrap draws for each permutation $\pi$.

Table \ref{tab:sim2} shows the results for various choices of $n$. By comparing the covariances  (\ref{eq:cov}) and (\ref{eq:covperm}) in Theorem \ref{thm:CLT}, tests permuting $L_n$ and $M_n$ will not have the correct level under $\mathcal{H}_{\mu(P)}$ even asymptotically, though the magnitude of the error is more alarming for $L_n$ than for $M_n$ in this simulation. A permutation test using $S_n$ on the other hand will be asymptotically correct under $\mathcal{H}_{\mu(P)}$ by Proposition \ref{prop:gauss} because it is equivalent to a Gaussian prepivoted test statistic; note the equality of the  rejection rates for $S_n$ before and after Gaussian prepivoting. The Gaussian prepivoted transformations of $L_n$ and $M_n$ restore asymptotic validity under $\mathcal{H}_{\mu(P)}$, shown as $n$ increases in the second set of columns. Comparing the second and third sets of columns, bootstrap prepivoting outperforms Gaussian prepivoting, resulting in tests that come much closer to attaining the nominal level. For $L_n$ in particular, applying bootstrap prepivoting after Gaussian prepivoting further improves the error rate in small samples relative to either (a) simply applying bootstrap prepivoting to $L_n$ itself; or (b) simply applying Gaussian prepivoting to $L_n$ without applying the bootstrap. The superiority of the tests in the fourth set of columns (bootstrap prepivoting after Gaussian prepivoting) for both $S_n$ and $L_n$ align with Theorem \ref{thm:indep} and the discussion in \S \ref{sec:bg}. For $M_n$, while the bootstrap prepivoted test statistics both outperform the asymptotically valid Gaussian prepivoted transform there is no noticeable difference between the two modes of bootstrap prepivoting in the third and fourth set of columns.
\begin{table}
\caption{Simulated Type I error rates for permutation tests under $\mathcal{H}_{\mu(P)}$ in the multivariate, two-sample setting. For each simulation $n_1 = 0.3 n$, and the nominal level of the tests is $\alpha = 0.05$. Within each set of columns, from left to right the test statistics are the pooled Hotelling $T$-squared, unpooled Hotelling $T$-squared, and max absolute $t$ statistic. For any given sample size $n$, the maximal estimated standard error for the difference between two estimated rejection probabilities is 0.0073.}

\label{tab:sim2}
\begin{center}
\begin{tabular}{l c c c c c c c c c c c c}
&\multicolumn{3}{c}{Original} &\multicolumn{3}{c}{Gaussian Prepivot} &\multicolumn{3}{c}{Bootstrap Prepivot}  &\multicolumn{3}{c}{Boot-after-Gauss}\\
$n$&$L_n$&$S_n$&$M_{n}$&$L_n$&$S_n$&$M_{n}$&$L_n$&$S_n$&$M_{n}$&$L_n$&$S_n$&$M_{n}$\\
$150$ & 0.820 & 0.368 & 0.330 & 0.323 & 0.368 & 0.310 &0.159 & 0.159 & 0.105 & 0.118 & 0.159 & 0.117\\
$500$ & 0.751 & 0.173 & 0.201 & 0.160 & 0.173 & 0.191 &0.096 & 0.112 & 0.076 & 0.084 & 0.112 & 0.079\\
$1000$ &0.704 & 0.107 & 0.139 & 0.101 & 0.107 & 0.132 & 0.067 &0.077 & 0.059 & 0.063 & 0.077 & 0.061\\
$2000$&0.692 & 0.082 & 0.103 & 0.080 & 0.082 & 0.095 & 0.056 & 0.061 & 0.056 & 0.054 & 0.061 & 0.057\\
$5000$&0.686 & 0.065 & 0.082 & 0.066 & 0.065 & 0.077 & 0.055& 0.055 & 0.053 & 0.053 & 0.055 & 0.052\\
\end{tabular}
\end{center}
\end{table}

\subsection{The $k$-sample problem: Analysis of variance, Tukey-Kramer, and an alternative robust test} \label{sec:anova}

We now consider univariate responses in the multi-group setting, taking $k=4$. Common test statistics for inference on $\mathcal{H}_{\mu(P)}$ with $k>2$ are the $F$-statistic arising from an analysis of variance computation $F_n$; and the Tukey-Kramer range statistic $T_n$,
\begin{align*}
F_n(Z_\pi) &= \frac{n\{\hat{\mu}_\pi C_2\text{diag}(b^{1/2})\} \{\hat{\mu}_\pi C_2\text{diag}(b^{1/2})\}^\T/(k-1)}{\hat{\nu}_\pi};\\
T_n(Z_\pi) &= \underset{1\leq \ell \leq k(k-1)/2}{\max}\;\; \frac{n^{1/2}|\hat{\mu}_{\pi}C_1|_\ell}{\{(C_1^\T\hat{\Lambda}_\pi C_1)_{\ell \ell}\}^{1/2}},
\end{align*} where $\text{diag}(b^{1/2})$ is a $k\times k$ diagonal matrix with $b_i^{1/2} = (n_i/n)^{1/2}$ in the $i$th diagonal. The conventional reference distributions make an assumption of equality of variances across the $k$ groups, an assumption not embedded within $\mathcal{H}_{\mu(P)}$. Hence the usual reference distributions do not generally satisfy Theorem \ref{thm:general} in the heteroskedastic case. Unlike in previous examples, Gaussian prepivoting for $F_n$ and $T_n$ yields distributions with neither closed-form distribution functions nor fast numerical approximation algorithms. While the transformation (\ref{eq:asym}) can be approximated through Monte-Carlo simulation as described in Remark \ref{rem:MC}, we instead proceed applying bootstrap prepivoting to $T_n$ and $F_n$ to restore asymptotic correctness.

In addition to these two statistics, we also consider the robust statistic suggested in \citet[Equation 3.2]{chu13}, shown to provide a robust permutation test under $H_{\mu(P)}$ in their Theorem 3.1. To express this statistic in the form (\ref{eq:T}),  let $\hat{D}_{\pi}$ be a $k\times k$ matrix with identical columns where $\hat{D}_{\pi ij} = (n_i/\hat{\Sigma}_{\pi i})/\{\sum_{\ell = 1}^k(n_\ell/\hat{\Sigma}_{\pi \ell})\}$ for all $j$, and let $\hat{C}_\pi = I_{k\times k} - \hat{D}_\pi$. Observe that $\hat{C}_\pi$ is a random matrix of column contrasts. Denoting their robust $k$- sample statistic as $W_n(Z_\pi)$, 
  \begin{align*}
W_n(Z_\pi) = \sum_{i=1}^k(n_i/n)(1/\hat{\Sigma}_{\pi \ell})\{(n^{1/2}\hat{\mu}_\pi\hat{C}_\pi)_i\}^2 = \sum_{i=1}^k\frac{n_i}{\hat{\Sigma}_{\pi i}}\left(\hat{\mu}_{\pi i} - \frac{\sum_{\ell=1}^kn_\ell\hat{\mu}_{\pi i}/\hat{\Sigma}_{\pi \ell}}{\sum_{\ell=1}^kn_\ell/\hat{\Sigma}_{\pi \ell}}\right)^2.
\end{align*} 
$W_n$ also satisfies Condition \ref{cond:T}. For this statistic, Gaussian prepivoting through (\ref{eq:asym}) returns the test statistic $G_{k-1}(W_n)$, where $G_{df}$ is the distribution function of a $\chi^2$ random variable with $df$ degrees of freedom. This is a monotone increasing transformation of $W_n$, such that the permutation distributions of the original test statistic and the Gaussian prepivoted test statistic yield identical $p$-values. Proposition \ref{prop:gauss} thus provides an alternative justification for the success of this test statistic: it is equivalent to a Gaussian prepivoted test statistic. We investigate permutation tests with both $W_n$ and its bootstrap prepivoted modification.

We draw 5000 data sets each of size $n$, with $n_1=0.1n$, $n_2=0.2n$, $n_3=0.3n$, $n_4=0.4n$. Let $\sigma_1=0.70, \sigma_2=0.55, \sigma_3 = 0.40$, and $\sigma_4=0.25$. For $i=1,\ldots,4$, we take $X_{ij}$ independent and identically distributed as $X_{ij}\sim \exp\{\mathcal{N}(0, \sigma_i^2)\} - \exp(\sigma^2_i/2)$ for $j=1,...,n_i$. For each data set, we conduct permutation tests using the original and bootstrap prepivoted simulations based upon 999 draws from the permutation distribution. When bootstrap prepivoting is used, we take 200 draws for each permutation $\pi$.

Table \ref{tab:sim3} shows the results for various choices of $n$. As a consequence of Theorem \ref{thm:CLT}, permutation tests based upon $F_n$ and $T_n$ will be asymptotically invalid under $\mathcal{H}_{\mu(P)}$. This is reflected in the first set of columns, where these tests reject well above the nominal level even at $n=2000$. The permutation test using $W_n$ will result in an asymptotically correct procedure by Proposition \ref{prop:gauss} because it is equivalent to its Gaussian prepivoted transformation. The second set of columns reflects Proposition \ref{prop:boot}: bootstrap prepivoting has restored first-order correctness to the permutation tests using $F_n$ and $T_n$ under $\mathcal{H}_{\mu(P)}$, and the nominal levels of these tests approach $\alpha=0.05$ as $n$ increases. We also see that the bootstrap prepivoted version of $W_n$ provides an improvement in the error in rejection rate relative to the tests using $W_n$ itself, particularly at $n=100$. This highlights the benefits of applying bootstrap prepivoting even with test statistics whose permutation tests are already asymptotically valid for testing equality of distributional parameters, as suggested by Theorem \ref{thm:indep}.
\begin{table}
\caption{Simulated Type I error rates for permutation tests under $\mathcal{H}_{\mu(P)}$ in the $k$-sample, univariate setting. We set $k=4$. In each simulation $n_1 = 0.1 n$, $n_2=0.2n$, $n_3=0.3n$, $n_4=0.4n$, and the nominal level of the tests is $\alpha = 0.05$. Within each set of columns, left to right the test statistics are the $F$-statistic from an analysis of variance, the Tukey-Kramer max statistic, and the statistic of \citet[Equation 3.2]{chu13}. For any given sample size $n$, the maximal estimated standard error for the difference between two estimated rejection probabilities is 0.0060.}\label{tab:sim3}
\begin{center}
\begin{tabular}{l c c c c c c}
&\multicolumn{3}{c}{Original}  &\multicolumn{3}{c}{Bootstrap Prepivot} \\
$n$&$F_n$&$T_n$&$W_{n}$&$F_n$&$T_n$&$W_{n}$\\
$100$ & 0.218 & 0.227 & 0.120 & 0.084 & 0.088&0.091\\
$250$ & 0.206 & 0.220 & 0.076 & 0.070 & 0.070 & 0.064\\
$1000$ &0.204 & 0.219 & 0.062 & 0.056 & 0.059 & 0.055 \\
$2000$&0.205 & 0.217 & 0.053 & 0.053 & 0.054 & 0.050\\
\end{tabular}
\end{center}
\end{table}

\subsection{Additional examples}
In \S A of the appendix, we consider multivariate analysis of variance with both $k>2$ and $d>2$. The results are qualitatively similar to those presented thus far. Wilk's Lambda, the Pillai-Bartlett Trace, the Lawley-Hotelling Trace, and Roy's Greatest Root are all of the form required in Condition \ref{cond:T}, but permutation tests using the statistics themselves can be invalid under $H_{\mu(P)}$. We also consider testing equality of medians in the two-sample case by way of bootstrap prepivoting. While this restores asymptotic validity, the higher-order considerations in \S \ref{sec:higher} do not apply, as asymptotic expansions for the bootstrap distributions do not exist \citep{hal89, sak00}.

\section{Discussion}
The robust permutation tests presented in this work are exact under the null of equality of distributions, while providing asymptotically valid inference for equality of parameters so long as asymptotically linear estimators exist. The conditions in this work do not encompass permutation tests based upon $k$-sample $U$-statistics such as Wilcoxon's rank sum test; see \citet{chu16ustat} for a discussion of robust permutation tests for two-sample $U$-statistics. There it is shown that while qualitatively similar ideas work for constructing robust permutation tests, the proofs vary considerably between the case of asymptotically linear estimators and the case of $U$-statistics.

The remedies presented in \S\S\ref{sec:asym}-\ref{sec:boot} provide two paths forward, one permuting analytical $p$-values based upon asymptotic normality and the other permuting bootstrap $p$-values. These are by no means the only available options. For instance, instead of the bootstrap, one could also construct $p$-values using subsampling. For the nonparametric Behrens-Fisher problem presented in \S \ref{sec:bf}, one could also permute the $p$-value from Welch's unpooled $t$-test using the Welch-Satterthwaite degrees of freedom. While not explicitly explored here, our use of permutation tests based upon asymptotically valid $p$-values has immediate consequences for power against contiguous alternatives to $\mathcal{H}_{\theta(P)}$: the prepivoted permutation tests have the same limiting local power as the large-sample tests whose $p$-values they permute. See \citet[Remark 2.3]{chu13} for further discussion. Hence there is nothing lost in terms of first-order power properties for the prepivoted permutation tests relative to large-sample alternatives, while much is gained by providing an exact test for $P_1=\ldots=P_k$.

The possibility of higher-order corrections described in \S\ref{sec:higher} provides theoretical motivation for using bootstrap $p$-values within permutation tests instead of other competitors: tests using bootstrap prepivoting maintain exactness under equality of distribution, while the order of approximation may be driven down under equality of parameters through bootstrap iteration. The computation required by embedding a bootstrap scheme within a permutation test makes repeated iterations of the bootstrap practically infeasible. When available, bootstrapping a statistic satisfying the conditions of Theorem \ref{thm:indep2} can provide an improvement in the order of approximation without requiring an additional layer of bootstrapping. An interesting area for future research is the potential use of weighted bootstrap iteration \citep{lee03} within permutation tests, which is known to provide subsequent improvements of $O(n^{-1})$ rather than $O(n^{-1/2})$ with each bootstrap iteration.

As previously mentioned, the primary downside of the method lies in the additional computation time required for enumerating the permutation distribution. Permuting analytical $p$-values as in Gaussian prepivoting requires computation on par with computing a bootstrap $p$-value, while permutating bootstrap $p$-values requires computation on par with double bootstrapping without providing the iterative improvement in the order of approximation that a $p$-value based upon a double bootstrap would often confer. The benefit provided by this additional layer of Monte-Carlo sampling is exactness under equality in distributions. Importantly, this property can \textit{only} be provided by permutation tests \citep[Theorem 5.8.1]{leh05}, such that the additional level of Monte-Carlo sampling is unavoidable in the setting considered here. A natural question, left for future research, is whether prepivoted permutation tests provide finite-sample improvements over the the test whose $p$-value it permutes when the distributions $P_1,...,P_k$ are ``similar'' in some sense yet unequal. 

\setcounter{table}{0}
\setcounter{equation}{0}
\newpage
\appendix
\renewcommand{\thetable}{A\arabic{table}}
\renewcommand{\theequation}{A\arabic{equation}}
\section{Additional simulation studies}
\subsection{Multivariate analysis of variance with prepivoted permutation tests} \label{sec:manova}

When $k>2$ and $d>2$, the usual tests for $\mathcal{H}_{\mu(P)}$ are based upon the eigenvalues of the product of the model sum of squares matrix and the inverse of the residual sum of squares matrix. Using the notation established in \S \ref{sec:notation},  this product may be written as 
\begin{align*} D_n(Z_\pi) &= n\left\{\hat{\mu}_\pi C_2\text{diag}(b^{1/2})\} \{\hat{\mu}_\pi C_2\text{diag}(b^{1/2})\}^\T/(k-1)\right\}\hat{\nu}_\pi^{-1}.
\end{align*}
Let $\lambda_{\pi j}$ be the $j$th eigenvalue of $D_n(Z_\pi)$. We consider the following common test statistics:
\begin{align*}
\text{Wilk's Lambda} &= \prod_{j=1}^k\{1/(1+\lambda_{\pi j})\};\\
\text{Pillai-Bartlett Trace} &= \sum_{j=1}^k\{\lambda_{\pi j}/(1+\lambda_{\pi j})\};\\
\text{Lawley-Hotelling Trace} &= \sum_{j=1}^k \lambda_{\pi j};\\
\text{Roy's Greatest Root} &= \underset{1\leq j \leq k}{\max}\;\; \lambda_{\pi j}.
\end{align*}

For each of these test statistics, we consider permutation inference using both the statistic itself and the bootstrap prepivoted transformations of these statistics. We note that unlike the other statistics considered in the main text, for Wilk's Lambda evidence against the null is suggested by small values of the test statistic; we simply take the negative of the test statistic and compute the right-tail permutation $p$-value as discussed in the manuscript. The other three tests reject with large values of the test statistic.

We simulate 5000 data sets of size $n$ with $k=4$ groups, with $n_1=0.1n$, $n_2=0.2n$, $n_3=0.3n$, $n_4=0.4n$. The outcome variable is of dimension $d_P=10$. Let $R(\rho)$ denote a $10\times 10$ correlation matrix with equal correlations $\rho$. Let $\rho_1 = 0.3, \rho_2 = 0.5, \rho_3 = 0.7, \rho_4 = 0.9$, and let $\sigma^2_1=1, \sigma^2_2=0.8, \sigma^2_3 = 0.6, \sigma^2_4 = 0.4$. For $j=1,\ldots,n_i$ and $i=1,\ldots,4$, we generate outcome variables as 
\begin{align*}
X_{ij} \sim \exp(\mathcal{N}_{10}\left(0, \sigma^2_iR(\rho_i)\right) - (\exp(\sigma^2_i/2),\ldots,\exp(\sigma^2_i/2))^\T,
\end{align*}
with the vector $(\exp(\sigma^2_i/2),\ldots,\exp(\sigma^2_i/2))^\T$ being of length 10. In each data set we conduct permutation inference using 999 permutations. When bootstrap prepivoting is applied, we use 200 bootstrap draws for each permutation. We investigate performance at $n=200, 1000$ and $2000$.

Table \ref{tab:sim4} shows the results. We see that the permutation tests using the original test statistics have rejection rates far above the nominal level even at $n=2000$, and by Theorem \ref{thm:CLT} this would persist asymptotically. The bootstrap prepivoted tests approach the nominal level as $n$ increases, reflecting Proposition \ref{prop:boot}.

\begin{table}
\caption{Simulated Type I error rates for permutation tests under $\mathcal{H}_{\mu(P)}$ in the $k$-sample, multivariate setting. We set $k=4$. In each simulation $n_1 = 0.1 n$, $n_2=0.2n$, $n_3=0.3n$, $n_4=0.4n$, and the nominal level of the tests is $\alpha = 0.05$. Within each set of four columns, left to right the test statistics are Wilk's Lambda, the Pillai-Bartlett Trace, the Lawley-Hotelling Trace, and Roy's Greatest Root. For any given sample size $n$, the maximal estimated standard error for the difference between two estimated rejection probabilities is 0.0069.}\label{tab:sim4}
\begin{center}
\begin{tabular}{l c c c c c c c c}
&\multicolumn{4}{c}{Original}  &\multicolumn{4}{c}{Bootstrap Prepivot} \\
$n$&W&P-B&L-H&R&W &P-B&L-H&R\\
$200$ & 0.731 & 0.721 & 0.741 & 0.774 & 0.107&0.101& 0.103 & 0.125\\
$1000$  & 0.723 & 0.722 & 0.726 & 0.781 & 0.074&0.068& 0.068 & 0.081\\
$2000$ & 0.708 & 0.708 & 0.709 & 0.783 & 0.066&0.062& 0.062 & 0.068\\
\end{tabular}
\end{center}
\end{table}


\subsection{Equality of medians} \label{sec:median}
When testing equality of medians in the two-sample setting, \citet[Example 2.2 and Section 4]{chu13} suggest permuting a studentized difference in medians using the bootstrap standard error. Here we demonstrate that the bootstrap could instead be used to prepivot the unstudentized difference in medians, providing an alternative to studentization for restoring asymptotic validity for permutation tests performed in the absence of group invariance. Let $\hat{m}_{\pi i}$ be the median of the responses in group $i$, and let $\hat{v}_{\pi i}$ be the bootstrap estimator of the variance for $n^{1/2}(\hat{m}_{\pi 1} - \hat{m}_{\pi 2})$; see \citet[\S 4]{chu13} for a closed form expression for $\hat{v}_\pi$. Define
\begin{align*}T_n(Z_\pi) &= {n^{1/2}(\hat{m}_{\pi 1} - \hat{m}_{\pi 2})};\\
S_n(Z_\pi) &= T_n(Z_\pi)/\hat{v}^{1/2}_{\pi i}.\\
\end{align*}We will compare inference based upon $T_n$ and $S_n$ to tests applying bootstrap prepivoting to $T_n$ and $S_n$.  We note that while the developments in \S\ref{sec:asym}-\S\ref{sec:boot} apply to both statistics, the discussion in \S\ref{sec:higher} does not apply to the unstudentized median or to the studentized median. While the bootstrap distributional estimates are consistent, subsequent terms in the asymptotic expansions do not exist \citep{hal89, sak00}. Bootstrap tests using the studentized and unstudentized statistics have the same order of error rejection probability as the large-sample test, and bootstrap iteration has absolutely no effect on the order of approximation. So while it is unclear whether the bootstrap will provide any benefit over studentization, it nonetheless will restore asymptotic validity by Proposition \ref{prop:boot}. 

We simulate 10,000 data sets. In each data set, we let $X_{1j}$ be independent and identically distributed draws from a standard normal distribution for $j=1,..,n_1$, and let $X_{2j}$ be independent and identically distributed draws from a normal distribution with mean zero and standard deviation 5 for $j=1,..,n_2$. For each data set permutation inference is based upon $999$ permutations. For bootstrap prepivoted test statistics, 500 bootstrap draws are used for each permutation.

Table \ref{tab:sim5} show the results for various choices of $n$. We first note that unlike the case of the difference in means, the permutation test using $T_n$ does not yield an asymptotically valid test for equality of medians even when $n_1=n_2$. We see that $S_n$ yields an asymptotically valid test by Proposition \ref{prop:gauss} because it is equivalent to a Gaussian prepivoted test statistic, and that the bootstrap prepivoted version of $T_n$ and $S_n$ are also asymptotically valid through an application of Proposition \ref{prop:boot}. 
\begin{table}
\caption{Simulated Type I error rates for permutation tests testing equality of medians in the two-sample, univariate case for various values of $n_1$ and $n_2$. The nominal level of the tests is $\alpha = 0.05$. For any given sample size $(n_1,n_2)$, the maximal estimated standard error for the difference between two estimated rejection probabilities is 0.0043.}\label{tab:sim5}
\begin{center}
\begin{tabular}{l c c c c c c}
&\multicolumn{2}{c}{Original}  &\multicolumn{2}{c}{Bootstrap Prepivot} \\
$(n_1,n_2)$&$T_n$&$S_n$&$T_n$&$S_n$\\
$(13, 13)$ & 0.223 & 0.123 & 0.139 & 0.099 \\
$(51, 51)$ & 0.119 & 0.062 & 0.072 & 0.050 \\
$(201, 201)$ &0.230 & 0.058 & 0.073 & 0.061 \\
$(401, 401)$ &0.223 & 0.053 & 0.064 & 0.059 \\
\end{tabular}
\end{center}
\end{table}


\newpage

\section{Pseudocode for bootstrap prepivoted permutation tests}
The following pseudocode implements a bootstrap prepivoted permutation test. To perform Gaussian prepivoting, simply replace $\hat{J}_{T_n}$ with $K_{T_n}$ as defined in (\ref{eq:asym}), which eliminates the need for any draws from the bootstrap distribution.

\begin{algorithm} 
	\SetAlgoLined
			\KwIn{The observed data $Z$; the test statistic $T_n(\cdot)$ and its bootstrap modification $\breve{T}_n(\cdot, \cdot)$; the number of bootstrap and permutation draws \texttt{nboot} and \texttt{nperm}.}
			\KwResult{The $p$-value for the bootstrap prepivoted test statistic}

			Compute $T_n(Z)$.
			
			\For{b=1 to b=\texttt{nboot}}{
				Draw $n_i$ observations $Z_{i1}^*,\ldots,Z_{in_i}^*$ from $\hat{P}_i$ (the empirical distribution of $Z_i$) for $i=1,..,k$.\\
				Stack $Z^*_{11},\ldots,Z^*_{kn_k}$ as a matrix $Z^*$.\\
				Compute $t_b = \breve{T}_n(Z^*, Z)$.}
			Compute 
			\begin{align*} \hat{J}_{T_n}\{T_n(Z), \hat{P}\} &= \frac{1}{\mathtt{nboot}}\sum_{b=1}^{\mathtt{nboot}}1\{t_b \leq T_n(Z)\}\end{align*}
			\For{p=1 to p=\texttt{nperm}}{
				Draw $\pi$ uniformly from the set of permutations $\mathcal{P}_n$.\\
				Create $Z_\pi$ by rearranging the rows of $Z$ based upon the permutation $\pi$.\\
				Compute $T_n(Z_\pi)$.\\
				\For{b=1 to b=\texttt{nboot}}{
				Draw $n_i$ observations $Z_{\pi i1}^*,\ldots,Z_{\pi in_i}^*$ from $\hat{P}_{\pi}$ (the empirical distribution of $Z_{\pi i}$) for $i=1,..,k$.\\
				Stack $Z^*_{\pi 11},\ldots,Z^*_{\pi kn_k}$ as a matrix $Z^*_\pi$.\\
				Compute $t_{\pi b} = \breve{T}_n(Z^*_\pi, Z_\pi)$.}
			Compute 
			\begin{align*} j_{p} &= \hat{J}_{T_n}\{T_n(Z_\pi), \hat{P}_\pi\} =  \frac{1}{\mathtt{nboot}}\sum_{b=1}^{\mathtt{nboot}}1\{t_{\pi b} \leq T_n(Z_\pi)\}\end{align*}
			}
			\Return{
			\begin{align*}
			p_{val} = \frac{1}{1+\mathtt{nperm}}\left(1 +  \sum_{p=1}^{\mathtt{nperm}}1[j_p \geq \hat{J}_{T_n}\{T_n(Z), \hat{P}\}]\right)
			\end{align*}	}
				
		\end{algorithm}
		\newpage 
\section{Bootstrap iteration}

Through iterated applications of bootstrap prepivoting one can attain further refinements to the order of approximation for permutation tests conducted when only $\mathcal{H}_{\theta(P)}$ holds if the expansions admit additional terms. Begin once again with a test statistic $T_n$ of the form (\ref{eq:T}). The first application of bootstrap prepivoting yields the test statistic $J_{T_{n}}\{T_{n}(Z_{\pi}), Z_{\pi}\}$, the transform of $T_n(Z_\pi)$ by its bootstrap distribution function. In the second round of bootstrap prepivoting, $J_{T_{n}}\{T_{n}(Z_{\pi}), Z_{\pi}\}$ \textit{itself} is transformed by a bootstrap estimate of its distribution function for each $\pi \in \mathcal{P}_{n}$. For a given permutation $\pi\in \cP_{n}$, one first constructs $\hat{P}^*_{\pi i}$ as the empirical distribution of $Z^*_{\pi i1},\ldots,Z^*_{\pi in_i}$ (the bootstrap samples from $\hat{P}_{\pi i}$) for $i=1,\ldots,k$. Then, for each $i$, one takes $n_i$ independent and identically distributed draws from $\hat{P}^*_{\pi i}$, call them $Z^{**}_{\pi i1},\ldots,Z^{**}_{\pi in_i}$, in so doing taking a bootstrap sample from the initial bootstrap sample in group $i$. Let $Z^{**}_{\pi i}$ be the $n_i \times {d}_P$ matrix whose $j$th row contains $Z^{**}_{\pi ij}$, and write $Z^{**}_{\pi}$ for the $n\times {d}_P$ matrix stacking $Z^{**}_{\pi 1}, ,\ldots,  Z^{**}_{\pi k}$ on top of one another. Define $J^{(1)}_{J_{n}}(x, \hat{P}_\pi)$ as
\begin{align*}
J^{(1)}_{J_n}(x, \hat{P}_\pi) = \text{pr}[J_{T_n}\{\breve{T}_n(Z_\pi^{**}, Z^*_\pi), \hat{P}^{*}_{\pi}\} \leq x \mid  Z_\pi].
\end{align*}   Regardless of whether or not Theorems \ref{thm:dep} or \ref{thm:indep} applied to the original test statistic, under Conditions \ref{cond:expand}-\ref{cond:dhat} and after the first application of bootstrap prepivoting we have that $H_{J_n}(x, P) = U(x) + O(n^{-j_1/2})$ and $R_{J_n}(x,  \hat{P}) = U(x) + O_P(n^{-j_1/2})$ for $j_1 = j+1$ under the setting of Theorem \ref{thm:indep} or $j_1=1$ under the setting of Theorem \ref{thm:dep}. Denote now $H^{(1)}_{J_n} = H_{J_n}$ and $R^{(1)}_{J_n} = R_{J_n}$ to reinforce that these are the true and permutation distributions after a single application of bootstrap prepivoting and that $J^{(1)}_{J_n}(x, P)$ and $J^{(1)}_{J_n}(x, P_\Pi)$ estimate these distributions. Supposing the existence of higher-order expansions for $H^{(1)}_{J_n}$, $R^{(1)}_{J_n}$ and their bootstrap analogues, in regular cases they would take the form
\begin{align*}
H^{(1)}_{J_{n}}(x, P) &= U(x) + n^{-j_1/2}\tilde{a}_{T_n}(x, P) + O(n^{-(j_1+1)/2}); \\
R^{(1)}_{J_{n}}(x,\hat{P}) &= U(x) + n^{-j_1/2}\tilde{b}_{T_n}(x,\hat{P}) + O_P(n^{-(j_1+1)/2});\\
J^{(1)}_{J_{n}}(x,\hat{P}) &= U(x) + n^{-j_1/2}\tilde{k}_{T_n}(x,\hat{P}) + O_P(n^{-(j_1+1)/2});\\
 J^{(1)}_{J_{n}}(x,\hat{P}_\Pi) &= U(x) + n^{-j_1/2}\tilde{k}_{T_n}(x,\hat{P}_\Pi) + O_P(n^{-(j_1+1)/2}).
\end{align*}
Consider the test statistic $J^{(1)}_{J_n}[J_{T_n}\{T_{n}(Z_\pi),\hat{P}_\pi\}, \hat{P}_\pi]$ and suppose that $\tilde{a}_{T_n}(x, P) - \tilde{k}_{T_n}(x,  \hat{P})$ and $\tilde{b}_{T_n}(x,\hat{P}) - \tilde{k}_{T_n}(x,\hat{P}_\Pi)$ are both $O_P(n^{-1/2})$ uniformly over $x$. The argument underpinning Theorem \ref{thm:indep} would then apply because the lead terms in the expansions all equal $U(x)$. Letting $H^{(2)}_{J_n}$ and $R^{(2)}_{J_n}$ be the true and permutation distributions for $J^{(1)}_{J_n}[J_{T_n}\{T_{n}(Z_\pi),\hat{P}_\pi\},\hat{P}_\pi]$, the statistic after two applications of bootstrap prepivoting, we would have $H^{(2)}_{J_n}(x, P) - R^{(2)}_{J_n}(x, \hat{P}) = O_P(n^{-(j_1+1)/2})$. This translates into an improvement from $O(n^{-j_1/2})$ to $O(n^{-(j_1+1)/2})$ in the error in rejection probability under $\mathcal{H}_{\theta(P)}$ when permuting $J^{(1)}_{J_n}[J_{T_n}\{T_{n}(Z_\pi),\hat{P}_\pi\}, \hat{P}_\pi]$.

Further iterations of bootstrap prepivoting may be applied, so long as the asymptotic expansions exist in the required form at each iteration. For instance, the third application of bootstrap prepivoting applies the bootstrap transform \begin{align*}
J^{(2)}_{J_n}(x, \hat{P}_\pi) = \text{pr}(J^{(1)}_{J_n}[J_{T_n}\{\breve{T}_n(Z^{***}_\pi, Z^{**}_\pi), \hat{P}^{**}_{\pi}\}, \hat{P}^*_\pi] \leq x \mid  Z_\pi),
\end{align*} where for $i=1,\ldots,k$ $Z^{***}_{\pi i 1},\ldots,Z^{***}_{\pi i n_i}$ are independent and identically distributed draws from $\hat{P}^{**}_{\pi i}$, the empirical distribution of $Z^{**}_{\pi i1}, \ldots, Z^{**}_{\pi in_i}$, and $\hat{P}^{**}_\pi = \prod_{i=1}^k\prod_{j=1}^{n_i}\hat{P}^{**}_{\pi i}$. In general, the $\ell$th application yields $H^{(\ell)}_{J_n}(x, P) - R^{(\ell)}_{J_n}(x, P) = O_P(n^{-(j_1+\ell-1)/2})$ uniformly in $x$ if the necessary expansions exist. 
\section{Theorem 1 and the role of contrast matrices}

Under Conditions 1 and 2, unconditional asymptotic normality for $V_n(Z) = \text{vec}\{\hat{\Theta}(Z)\hat{C}(Z)\}$ follows immediately from the central limit theorem and Slutsky's theorem. The form of the covariance given in (\ref{eq:cov}) may be derived using \citet[Equation 6]{hen79}. To study the probability limit of $R_{V_n}(x)$ at all points $x$, we use a technique devised by \citet{hoe52} and generalized by \citet{dum13} and  \citet{chu13, chu16} which characterizes weak convergence in probability of permutation distributions in terms of the joint unconditional convergence of suitable test statistics. Adapting Lemma A.1 of \citet{chu16} to our setting: 
\begin{lemma}\label{lem:hoeffding}
Let $\Pi$ and $\Pi'$ be independent and identically distributed uniformly over $\mathcal{P}_n$, let $Z$ be distributed according to $P^n$ as before, and let $Z$ be independent of $\Pi$ and $\Pi'$. Then, $R_{V_n}$ converges weakly in probability to some law $Q$ if and only if $\{V_n(Z_\Pi), V_n(Z_{\Pi'})\}$ converges in distribution to $(U, U')$, where $U$ and $U'$ are independent and identically distributed with law $Q$.
\end{lemma}

Hence in order for $R_{V_n}$ to converge weakly in probability to a fixed law $Q$, $\{V_n(Z_\Pi), V_n(Z_{\Pi'})\}$ must converge jointly in distribution to random variables distributed according to $Q$ and must be asymptotically independent of one another. Under Conditions 1 and 2 the limit law $Q$ for $V_n$ must be that of a multivariate Gaussian. Therefore, asymptotic independence of the limiting random variables $U$ and $U'$ can be established by showing that $Cov(U_i, U'_{i'}) = 0$ for all $i,i'$.

Let $\Theta(\bar{P})$ be the $d\times k$ matrix with $\theta(\bar{P})$ in each of the $k$ columns, and let $\bar{Z}$ contain $n$ independent and identically distributed draws from the mixture distribution $\bar{P}$. We now apply the Cram\'er-Wold device to directly leverage the proof of Theorem 3.1 in \citet{chu13} for the case $d=1$. \citet{chu13} use slightly different scaling in their proof. To relate their results to ours, define $\tilde{\theta}_i(\bar{P}) = n_i^{1/2}\theta(\bar{P})$ and $\hat{\tilde{\theta}}(\bar{Z}_{\Pi i}) = n_i^{1/2}\hat{\theta}(\bar{Z}_{\Pi i})$. Let $\tilde{\Theta}(\bar{P})$ be the $d\times k$ matrix with $\tilde{\theta}_i(\bar{P})$ in the $i$th column, and define $\hat{\tilde{\Theta}}(Z_\Pi)$ analogously. 

The proof of Theorem 3.1 in \citet{chu13} establishes asymptotic normality for $[\text{vec}\{\hat{\tilde{\Theta}}(\bar{Z}_\Pi) - \tilde{\Theta}(\bar{P})\}, \text{vec}\{\hat{\tilde{\Theta}}(\bar{Z}_{\Pi'}) - \tilde{\Theta}(\bar{P})\}]^\T$ when $d=1$ under Conditions 1 and 2. Applying the Cram\'er-Wold device and recalling the relationship between $\hat{\Theta}(Z_\Pi)$ and $\hat{\tilde{\Theta}}(Z_\pi)$, it follows that $[n^{1/2}\text{vec}\{\hat{\Theta}(\bar{Z}_\Pi) - \Theta(\bar{P})\},n^{1/2}\text{vec}\{\hat{\Theta}(\bar{Z}_{\Pi'}) - \Theta(\bar{P})\}]^\T$ converges in distribution to a mean zero multivariate Gaussian $(\text{vec}(D), \text{vec}(D'))$, where $D$ and $D'$ are the limit distributions of $n^{1/2}\{\hat{\Theta}(\bar{Z}_\Pi) - \Theta(\bar{P})\}$ and $n^{1/2}\{\hat{\Theta}(\bar{Z}_{\Pi'}) - \Theta(\bar{P})\}$ respectively and are of dimension $d\times k$. The covariance matrices for $\text{vec}(D)$ and $\text{vec}(D')$ are identical and equal $\bar{\Gamma}$, a $kd\times kd$ block-diagonal matrix with $\bar{\Sigma}/p_i$ in the $i$th of $k$ blocks as described in the statement of Theorem \ref{thm:CLT}. On the off-diagonal, modifying (S25) of \citet{chu13} for our scaling we have for any columns $i, i'$ and any rows $j, j'$
\begin{align}
\text{cov}(D_{j i}, D'_{j' i'}) &= \bar{\Sigma}_{j j'}, \label{eq:covsup}
\end{align} where $D_{j i}$ and $D'_{j' i'}$ are the $\{j,i\}$ and $\{j',i'\}$ entries of $D$ and $D'$. 

Lemma \ref{lem:hoeffding} cannot be be applied directly because (\ref{eq:covsup}) does not equal zero. Recall from Condition \ref{cond:T} that $\bar{C}$ is a matrix of column contrasts and is assumed be the probability limit of $\hat{C}(Z_\Pi)$. We now consider the limit distribution of $(n^{1/2}\text{vec}[\{\hat{\Theta}(Z_\Pi) - \Theta(\bar{P})\}\bar{C}], n^{1/2}\text{vec}[\{\hat{\Theta}(\bar{Z}_{\Pi'}) - \Theta(\bar{P})\}\bar{C}])^\T$, which equals $[n^{1/2}\text{vec}\{\hat{\Theta}(\bar{Z}_\Pi)\bar{C}\}, n^{1/2}\text{vec}\{\hat{\Theta}(\bar{Z}_{\Pi'})\bar{C}\}]^\T$ because $\bar{C}$ is a matrix of column contrasts and the columns of $\Theta(\bar{P})$ are identical. This vector converges in distribution to $\{\text{vec}(D\bar{C}), \text{vec}(D'\bar{C})\}^\T$ with $(D,D')$ defined as before. Choose any two columns $\ell$ and $\ell '$ of $\bar{C}$ and any two rows $j$ and $j'$ of $\hat{\Theta}(\bar{Z}_\pi)$, and consider the limiting value of the covariance between $n^{1/2}\sum_{i=1}^k\hat{\theta}_j(\bar{Z}_{\Pi i})\bar{C}_{i \ell}$ and $n^{1/2}\sum_{i=1}^k\hat{\theta}_{j'}(\bar{Z}_{\Pi' i})\bar{C}_{i \ell' }$, i.e. $\text{cov}(\sum_{i=1}^kD_{j i}\bar{C}_{i \ell }, \sum_{i = 1}^k D'_{j' i}\bar{C}_{ i\ell'})$.
\begin{align}
\text{cov}\left\{n^{1/2}\sum_{i=1}^k\hat{\theta}_j(\bar{Z}_{\Pi i})\bar{C}_{i\ell}, n^{1/2}\sum_{i=1}^k\hat{\theta}_{j'}(\bar{Z}_{\Pi' i})\bar{C}_{i \ell'}\right\} &\rightarrow \bar{\Sigma}_{jj'}\sum_{i=1}^k\bar{C}_{i\ell}\sum_{i=1}^k \bar{C}_{i\ell'}\label{eq:contrast}\\
&= 0,\nonumber
\end{align}
where the last line follows because $\bar{C}$ is a matrix of column contrasts. This holds for any columns $\ell$ and $\ell'$ of $\bar{C}$ and for any $\hat{\theta}_j(\bar{Z}_\Pi)$, $\hat{\theta}_{j'}(\bar{Z}_{\Pi'})$. As a result, $\text{cov}[n^{1/2}\text{vec}\{\hat{\Theta}(\bar{Z}_\Pi)\bar{C}\}, n^{1/2}\text{vec}\{\hat{\Theta}(\bar{Z}_{\Pi'})\bar{C}\}]$ converges to a $kd\times kd$ matrix containing all zeroes. Using asymptotic normality, it follows that $\text{vec}\{n^{1/2}\hat{\Theta}(\bar{Z}_\Pi)\bar{C}\}$ and $\text{vec}\{n^{1/2}\hat{\Theta}(\bar{Z}_{\Pi'})\bar{C}\}$ are asymptotically independent. The form of the covariance given in (\ref{eq:covperm}) may again be derived from \citet[Equation 6]{hen79}.

The above proof applied to independent permutations of $\bar{Z}$, containing $n$ independent and identically distributed draws from $\bar{P}$. To complete the proof for $Z$, containing $n_i$ samples from $P_i$ for $i=1,\ldots,k$, we appeal to the coupling construction described in Section 5.3 of \citet{chu13}; see the proof of Theorem 3.1 of \citet{chu13} for details. Lemma \ref{lem:hoeffding} along with Slutsky's theorem for permutation distributions may then be applied to conclude that that the permutation distribution of $R_{V_n}$ converges weakly in probability to the law of a mean zero multivariate Gaussian with covariance (\ref{eq:covperm}) as desired. 

The above derivation shows the necessity of the columns of $\bar{C}$ being contrast matrices under our conditions. Suppose some column $\ell$ of $\bar{C}$ was not a contrast. Then, (\ref{eq:contrast}) with $\ell=\ell'$ and taking $j=j'$ for any $j$ would only be guaranteed to be zero in the degenerate case where $\bar{\Sigma}_{jj}=0$ for all $j$, which is disallowed under Condition 2. If (\ref{eq:contrast}) does not equal zero, asymptotic independence cannot be achieved and hence the permutation distribution for $n^{1/2}\text{vec}\{\hat{\Theta}(Z)\hat{C}(Z)\}$ cannot converge weakly in probability to the law of any fixed distribution by Lemma \ref{lem:hoeffding}. 
 
\section{Theorem 2}
We prove the statement for $H_{F_{n}}(x, P)$, the proof for $R_{F_n}(x, P)$ being analogous. Let $F_{T_n}^{-1}(x, \hat{P}) = \inf\{x : F_{T_n}(x, \hat{P})\geq x\}$ be defined as before, and define $\tilde{F}_{T_n}^{-1}(x, \hat{P}) = \sup\{x : F_{T_n}(x, \hat{P})\leq x\}$. Recalling that $H_{F_{n}}(x, P) = \P[F_{T_{n}}\{T_{n}(Z),  \hat{P}\} \leq x]$, we have
\begin{align*}
\P\{T_{n}(Z) \leq \tilde{F}_{T_{n}}^{-1}(x,\hat{P})\}\leq \P[F_{T_{n}}\{T_{n}(Z),  \hat{P}\} \leq x]\leq \P\{T_{n}(Z) \leq F_{T_{n}}^{-1}(x,\hat{P})\}.
\end{align*}
As $H_T(x,P)$ is continuous and strictly increasing as a function of $x$, we have by Lemma 11.2.1 of \citet{leh05} that $F_{T_{n}}^{-1}(x,\hat{P})$ converges in probability to $H_T^{-1}(x, P)$ for any $x\in [0,1]$. Analogous steps show that $\tilde{F}_{T_{n}}^{-1}(x,\hat{P})$ also converges in probability to $H_T^{-1}(x, P)$ for any $x\in [0,1]$.  Corollary 11.2.3 of \citet{leh05} then yields that $\P\{T_{n}(Z) \leq F_{T_{n}}^{-1}(x,\hat{P})\}$ and $\P\{T_{n}(Z) \leq \tilde{F}_{T_{n}}^{-1}(x,\hat{P})\}$ both converge to $H_T\{H_T^{-1}(x,P), P\} = U(x)$, the distribution function of a uniform random variable on [0,1], implying convergence of $\P[F_{T_{n}}\{T_{n}(Z),  \hat{P}\} \leq x]$ to $U(x)$ at all points $x\in[0,1]$. Applying Polya's theorem \citep[Theorem 11.2.9]{leh05} completes the proof.

\section{Proof of Proposition 1}
Under Condition \ref{cond:T}, the function
\begin{align*}h(\Gamma_0, n^{1/2}\Theta_0C_0, \eta_0, C_0) &= \gamma^{(kd)}_{0, \Gamma_0}\left\{{a}: {g}(\text{vec}_{d,k}^{-1}(a)C_0, \eta_0) \leq g(n^{1/2}\Theta_0C_0, \eta_0)\} \right\} 
\end{align*}
is jointly continuous in $\Gamma_0, n^{1/2}\Theta_0C_0, \eta_0, C_0$; see Lemma B of \citet{coh21gauss} for a proof. Moreover, observe that by Theorem \ref{thm:CLT} and the continuous mapping theorem, $h\{\Gamma, n^{1/2}\hat{\Theta}(Z)C, \eta, C\}$ converges in distribution to a standard uniform, and the permutation distribution of $h\{\bar{\Gamma}, n^{1/2}\hat{\Theta}(Z)\bar{C}, \bar{\eta}, \bar{C}\}$ converges weakly in probability to the law of a standard uniform distribution. Note that $K_{T_n}\{T_n(Z_\pi), \hat{P}_\pi\} = h\{\hat{\Gamma}(Z_\pi), n^{1/2}\hat{\Theta}(Z_\pi)\hat{C}(Z_\pi), \hat{\eta}(Z_\pi), \hat{C}(Z_\pi)\}$. That  $\{\hat{C}(Z), \hat{\eta}(Z), \hat{\Gamma}(Z)\}$ converges in probability to $(C, \eta, \Gamma)$ by  Conditions \ref{cond:T} and \ref{cond:sigma} provides the conclusion of the proposition pertaining to the unconditional law $H_{K_n}$ through Theorem 1, Slutsky's theorem, and continuous mapping theorem. Next, $\{\hat{C}(Z_\Pi), \hat{\eta}(Z_\Pi), \hat{\Gamma}(Z_\Pi)\}$ converges in probability to $(\bar{C}, \bar{\eta}, \bar{\Gamma})$ by  Conditions \ref{cond:T} and \ref{cond:sigma}, with the result for $\hat{\Gamma}(Z_\Pi)$ using contiguity results in Lemma 5.3 of \citet{chu13}, available to us under Condition \ref{cond:frac}. Hence, Slutsky's theorem and the continuous mapping theorem for permutation distributions \citep[Lemmas A.4 and A.5]{chu16} in concert with Theorem 1 then provide the result for the permutation distribution $R_{K_n}$. 

\section{Proof of Proposition 2}
\setcounter{equation}{2}

The proof closely follows that of Theorem 2.6 of \citet{chu16}, differing primarily in the condition used to ensure bootstrap consistency.  For each $i$ let $G_{n,i}(x, \hat{P}_{\pi i})$ be the bootstrap estimator for the distribution of $n^{1/2}\{\hat{\theta}(Z_{\pi i}) - \theta(\bar{P})\}$, i.e. the distribution of $n^{1/2}\{\hat{\theta}(Z_{\pi i}^*) - \hat{\theta}(Z_{\pi i})\}$ given $Z_{\pi}$ when $Z_{\pi i}^*$ contains $n_i$ independent and identically distributed draws from $\hat{P}_{\pi i}$. Fix a $\delta > 0$ and an $\varepsilon > 0$ and define the set $\mathcal{G}_n$ as
\begin{align*}
\mathcal{G}_n &= \{\pi: \underset{x}{\sup}\; |G_{n, i}(x, \hat{P}_{\pi  i}) - G_i(x, \bar{P})| < \delta,\;\;\; i = 1,\ldots,k\}\\& \cap \{\pi: \P\{|\hat{\eta}(Z^*_\pi) - \bar{\eta}| > \varepsilon  \wedge |\hat{C}(Z^*_\pi) - \bar{C}| > \varepsilon \mid Z_\pi\} < \delta\},
\end{align*}  
where $G_i(x, \bar{P})$ is the distribution function of a multivariate Gaussian with mean zero and covariance $\bar{\Sigma}/p_i$. Let $\mathcal{G}_n^c$ be the complement of $\mathcal{G}_n$. 

The randomization distribution $R_{J_n}(x, \hat{P})$ of $J_{T_n}\{T_n(Z_\pi), \hat{P}_\pi\}$ can be expressed as
\begin{align*}
R_{J_n}(x, \hat{P}) = \frac{1}{|\mathcal{P}_n|}\sum_{\pi \in \mathcal{G}_n} \1[J_{T_n}\{T_n(Z_\pi), \hat{P}_\pi\} \leq x] + \frac{1}{|\mathcal{P}_n|}\sum_{\pi \in \mathcal{G}^c_n} \1[J_{T_n}\{T_n(Z_\pi), \hat{P}_\pi\}\leq x].
\end{align*}
We begin by showing that $|\mathcal{G}_n|/|\mathcal{P}_n|$ converges to 1 in probability. To do so, we can show that for each $i$ and for any $\delta, \varepsilon > 0$,
\begin{align*}
\frac{1}{|\mathcal{P}_n|}\sum_{\pi \in \mathcal{P}_n}\1\{\underset{x}{\sup} |G_{n, i}(x, \hat{P}_{\pi i}) - G_i(x, \bar{P})| < \delta\} \overset{P}{\rightarrow} 1,
\end{align*}
and that
\begin{align*}
\frac{1}{|\mathcal{P}_n|}\sum_{\pi \in \mathcal{P}_n}1[\P\{|\hat{\eta}(Z^*_\pi) - \bar{\eta}| > \varepsilon  \wedge |\hat{C}(Z^*_\pi) - \bar{C}| > \varepsilon \mid Z_\pi\} < \delta] \overset{P}{\rightarrow} 1.
\end{align*}
It is sufficient to show that for each $i$
\begin{align} \P(\underset{x}{\sup}\; |G_{n, i}(x, \hat{P}_{\Pi i}) - G_i(x, \bar{P})| > \delta) \rightarrow 1, \label{eq:convergeG} \end{align} and that 
\begin{align}
\P[\P\{|\hat{\eta}(Z^*_\Pi) - \bar{\eta}| > \varepsilon  \wedge |\hat{C}(Z^*_\Pi) - \bar{C}| > \varepsilon \mid Z_\Pi\} < \delta] \rightarrow 1 \label{eq:convergeProb}
\end{align}
with $\Pi$ uniform over $\mathcal{P}_n$. (\ref{eq:convergeProb}) holds for any $\varepsilon, \delta$ by Condition \ref{cond:Tboot}. To prove (\ref{eq:convergeG}), we use the contiguity results presented in Lemma 5.3 of \citet{chu13} which are at our disposal under Condition \ref{cond:frac}. Letting $\hat{Q}_i$ be the empirical distribution of $Y_1,\ldots,Y_{n_i}$ when $Y_i$ are independent and identically distributed according to $\bar{P}$, Lemma 5.3 of \citet{chu13} shows that  (\ref{eq:convergeG}) holds if
\begin{align*} \P(\underset{x}{\sup}\; |G_{n, i}(x, \hat{Q}_{i}) - G_i(x, \bar{P})| < \delta) \rightarrow 1 \end{align*}
Under Conditions \ref{cond:asymlin} and \ref{cond:asymlinboot}, this is assured for each $i=1,\ldots,k$ by the bootstrap central limit theorem of \citet{liu89}. Hence, (\ref{eq:convergeG}) holds for any $\delta > 0$, implying the weak consistency of the bootstrap for $G_{n,i}(x, \hat{P}_{\Pi i})$ for all $i$. Furthermore, observe that given any $Z_\pi$, the random variables $n^{1/2}\{\hat{\theta}(Z_{\pi i}^*) - \hat{\theta}(Z_{\pi i})\}$ are jointly independent for $i=1,\ldots,k$ by construction. Therefore, an application of the Cram\'er-Wold device yields the joint convergence of $[n^{1/2}\{\hat{\theta}(Z_{\pi 1}^*) - \hat{\theta}(Z_{\pi 1})\},\ldots,n^{1/2}\{\hat{\theta}(Z_{\pi k}^*) - \hat{\theta}(Z_{\pi k})\}]^\T$ to $k$ independent multivariate Gaussians, with component-wise distribution functions $G_i(x, \bar{P})$ for $i=1,\ldots,k$.

For any $\delta>0$ and $\varepsilon > 0$ we can thus re-express $R_{J_n}$ as 
\begin{align}\label{eq:Rnew}
R_{J_n}(x, \hat{P}) = \frac{1}{|\mathcal{P}_n|}\sum_{\pi \in \mathcal{G}_n} \1[J_{T_n}\{T_n(Z_\pi), \hat{P}_\pi\} \leq x] + o_P(1).
\end{align}

Recall the construction of $\breve{T}_n(Z_\pi, Z^*_\pi)$ in (\ref{eq:Tboot}) used in defining $J_{T_n}$, 
\begin{align*}
\breve{T}_n(Z^*_\pi, Z_\pi) &= g\{n^{1/2}(\hat{\Theta}^*_\pi - \hat{\Theta}_\pi)\hat{C}^*_\pi, \hat{\eta}^*_\pi\},
\end{align*}
where the $i$th column of $\hat{\Theta}^*_\pi - \hat{\Theta}_\pi$ is $\hat{\theta}(Z_{\pi i}^*) - \hat{\theta}(Z_{\pi i})$.  Applying the continuous mapping theorem and Slutsky's theorem, we then have that $J_{T_n}(x, \hat{P}_\pi)$ converges in probability to $R_T(x, P)$ at all points $x$, where $R_T(x,P)$ is the distribution of $g(U, \bar{\eta})$ and $U$ has a multivariate normal distribution with covariance given by (\ref{eq:covperm}) in the manuscript.  Therefore, with probability tending to one, the remaining term in (\ref{eq:Rnew}) is bounded for any $\epsilon > 0$ as
\begin{align*}
&\frac{1}{|\mathcal{P}_n|}\sum_{\pi \in \mathcal{G}_n} \1[R_T\{T_n(Z_\pi), P\} \leq x - \epsilon ]\\\leq & \frac{1}{|\mathcal{P}_n|}\sum_{\pi \in \mathcal{G}_n} \1[J_{T_n}\{T_n(Z_\pi), \hat{P}_\pi\} \leq x]\\ \leq &\frac{1}{|\mathcal{P}_n|}\sum_{\pi \in \mathcal{G}_n} \1[R_T\{T_n(Z_\pi), P\} \leq x + \epsilon ].
\end{align*}

Using Theorem \ref{thm:CLT} and the continuous mapping theorem for permutation distributions \citep[Lemma A.6]{chu16}, we know that the permutation distribution of $T_n$, $R_{T_n}(x,\hat{P})$, also converges in probability to $R_T(x,P)$ at all points $x$, with $R_T(\cdot, P)$ continuous and strictly increasing at $R_T^{-1}(\cdot, P)$ by assumption. Again applying the continuous mapping theorem for permutation distributions, we have
\begin{align*}
\frac{1}{|\mathcal{P}_n|}\sum_{\pi \in \mathcal{G}_n} \1[R_T\{T_n(Z_\pi), P\} \leq x - \epsilon ] & \overset{P}{\rightarrow} x - \epsilon;\\
\frac{1}{|\mathcal{P}_n|}\sum_{\pi \in \mathcal{G}_n} \1[R_T\{T_n(Z_\pi), P\} \leq x + \epsilon ]  & \overset{P}{\rightarrow} x + \epsilon,\\
\end{align*}
such that for any $\epsilon > 0$ and with probability tending to one
\begin{align*}
x-\epsilon\leq  \frac{1}{|\mathcal{P}_n|}\sum_{\pi \in \mathcal{G}_n} \1[J_{T_n}\{T_n(Z_\pi), \hat{P}_\pi\} \leq x]\leq x+\epsilon.
\end{align*}
Sending $\epsilon$ to zero completes the proof.

\section{Proof of Theorem 3}
\begin{proof}
Observe that $K_{T_n}(x,  \hat{P}) = A_{T_n}(x, P) + O_P(n^{-1/2})$ and that $A_{T_n}(x, P) = H_{T_n}(x, P) + O_P(n^{-j/2})$. By continuity of $a_{T_n}(x, P)$ and (\ref{eq:H}), we have   $\P[H_{T_{n}}\{T_{n}(Z), P\}\leq x] = U(x) + O_P(n^{-j/2})$ uniformly in $x$. Therefore, we have
\begin{align*} \P[K_{T_n}\{T_{n}(Z),  \hat{P}\} \leq x] &= \P[H_{T_{n}}\{T_{n}(Z), P\} + O_P(n^{-1/2}) +  O_P(n^{-j/2}) \leq x]\\
&= U(x) + O(n^{-1/2})\end{align*}
uniformly in $x$. As a result,

\begin{align*} \P[J_{T_n}\{T_{n}(Z),  \hat{P}\} \leq x] &= \P[K_{T_{n}}\{T_{n}(Z), P\} + O_P(n^{-j/2}) \leq x]\\
&= U(x) + O(n^{-1/2})\end{align*}uniformly in $x$. The analogous derivation yields that
\begin{align*}
R_{J_n}(x, \hat{P}) &= \P\left[J_{T_{n}}\{T_{n}(Z_\Pi),\hat{P}_\Pi\} \leq x \mid Z\right]
\\&=  U(x) + O_P(n^{-1/2})
\end{align*}
uniformly in $x$, such that $H_{J_n}(x, P) - R_{J_n}(x,  \hat{P}) = O_P(n^{-1/2})$ uniformly in $x$.
\end{proof}

\section{Proof of Theorem 4}

\begin{proof}
As $A_{T_n}(x,P)=K_{T_n}(x,\hat{P})$, we have
\begin{align*}
H_{T_{n}}(x, P) - J_{T_{n}}(x,\hat{P}) = n^{-j/2}\{a_{T_n}(x,P) - k_{T_n}(x,\hat{P})\} + O_P(n^{-(j+1)/2}) = O_P(n^{-(j+1)/2})
\end{align*}
uniformly in $x$. Further observe that under the null, by continuity of $a_{T_n}(x,P)$ 
\begin{align*}
\P[H_{T_{n}}\{T_{n}(Z), P\} \leq x] = U(x) + O(n^{-(j+1)/2})
\end{align*}
uniformly in $x$. As a consequence,
  \begin{align*}
H_{J_n}(x, P) &= \P[J_{T_{n}}\{T_{n}(Z),  \hat{P}\} \leq x]\\&=\P[H_{T_{n}}\{T_{n}(Z),P\} \leq x + O_P(n^{-(j+1)/2})] = U(x) + O(n^{-(j+1)/2}),
\end{align*}
uniformly in $x$. An analogous derivation gives
\begin{align*}
R_{T_{n}}(x,  \hat{P}) - J_{T_{n}}(x,\hat{P}_\Pi)  = O_P(n^{-(j+1)/2});\;\;\;
\P[R_{T_{n}}\{T_{n}(Z),  \hat{P}\} \leq x\mid Z] = U(x) + O_P(n^{-(j+1)/2}),
\end{align*}
uniformly in $x$ such that
\begin{align*}
R_{J_n}(x,  \hat{P}) &= \P\left[J_{T_{n}}\{T_{n}(Z_\Pi),\hat{P}_\Pi\} \leq x \mid Z\right]
\\&= \P\left[R_{T_{n}}\{T_{n}(Z_\Pi), \hat{P}\} \leq x + O_P(n^{-(j+1)/2})\mid Z\right] = U(x) + O_P(n^{-(j+1)/2}).
\end{align*}
uniformly in $x$. This yields $H_{J_n}(x,P) - R_{J_n}(x,\hat{P}) = O_P(n^{-(j+1)/2})$ uniformly in $x$. 
\end{proof}

\bibliographystyle{apalike}
\bibliography{bibliography.bib}

\end{document}